\documentclass[12pt]{amsart}
\DeclareFontFamily{OML}{script}{}
\DeclareFontShape{OML}{script}{m}{it}
{ <5-20> rsfs10 }{}
\DeclareMathAlphabet{\mathscript}{OML}{script}{m}{it}

\renewcommand{\mathcal}[1]{{\mathscript #1}\hspace{0.2ex}}
\usepackage{calligra}
\usepackage{pstricks}
\usepackage{color}
\ifx\red\undefined
\newcommand{\red}{\color{red}}
\fi

\usepackage{graphicx,graphics}
\usepackage{cite}
\usepackage{pifont}
\usepackage{amsthm}

\usepackage{amscd}
\usepackage{amsmath}
\usepackage{latexsym}
\usepackage{amsfonts}
\usepackage{amssymb}
\usepackage{color}

\ifx\text\undefined
\newcommand{\text}{\mbox}
\fi
\ifx\operatorname\undefined
\newcommand{\operatorname}{\mathop}
\fi
\newcommand\be{\begin{equation}}
\newcommand\ee{\end{equation}}
\newcommand\bea{\begin{eqnarray}}
\newcommand\eea{\end{eqnarray}}
\newcommand\beaa{\begin{eqnarray*}}
\newcommand\eeaa{\end{eqnarray*}}

\newenvironment{eqa}{\begin{equation}%
  \begin{array}{rcl}}{\end{array}\end{equation}}
\newcommand\beqa{\begin{eqa}}
\newcommand\eeqa{\end{eqa}}
\numberwithin{equation}{section}

\newtheorem{thm}{Theorem}[section]

\newtheorem{lem}{Lemma}[section]

\newtheorem{rem}{Remark}[section]

\newcommand{\void}[1]{}

\numberwithin{equation}{section}

\begin{document}

\title[Asymptotic stability]{Asymptotic stability for a free boundary tumor model with angiogenesis}
\author[Huang]{Yaodan Huang}%黄瑶丹}
\author[Zhang]{Zhengce Zhang}
\author[Hu]{Bei Hu}

\date{\today}
\address[1]{School of Mathematics, Sun Yat-Sen University, Guangzhou, 510275, China; School of Mathematics and Statistics, Xi'an Jiaotong University,
Xi'an, 710049, China}
\address[2]{School of Mathematics and Statistics, Xi'an Jiaotong University,
Xi'an, 710049, China}
\address[3]{Department of Applied Computational Mathematics and Statistics, University of Notre Dame, Notre Dame, IN 46556, USA}
%\address[3]{Institute of Mathematical Sciences, Renmin University, Beijing 100872, China}
\email{huangyd35@mail.sysu.edu.cn, zhangzc@mail.xjtu.edu.cn, b1hu@nd.edu}
%\thanks{Corresponding author: Yaodan Huang}
%\thanks{This work was supported by the National Natural Science Foundation of China (No. 11371286, 11401458). }

\maketitle

\begin{abstract}
%A free boundary tumor model with vasculature which supplies nutrients to the tumor is considered. The tumor cells proliferate at a rate $\mu$. It was recently established in \cite{HZH2} that there exists a threshold value $\mu^\ast=\mu^\ast(R_S)$ such that the unique radially symmetric stationary solution is linearly stable under non-radial perturbations for $0<\mu<\mu^\ast$ and linearly unstable for $\mu>\mu^\ast$. In this paper, we show that for $0<\mu<\mu^\ast$, this stationary solution is actually asymptotically stable under nonlinear settings.
In this paper, we study a free boundary problem modeling solid tumor growth with vasculature which supplies nutrients to the tumor; this is characterized in the Robin boundary condition. It was recently established [Discrete Cont. Dyn. Syst. 39 (2019) 2473-2510] that for this model, there exists a threshold value $\mu^\ast$ such that the unique radially symmetric stationary solution is linearly stable under non-radial perturbations for $0<\mu<\mu^\ast$ and linearly unstable for $\mu>\mu^\ast$. In this paper we further study the nonlinear stability of the radially symmetric stationary solution, which introduces a significant mathematical difficulty: the center of the limiting sphere is not known in advance owing to the perturbation of mode 1 terms. We prove a new fixed point theorem to solve this problem, and finally obtain that the radially symmetric stationary solution is nonlinearly stable for $0<\mu<\mu^\ast$ when neglecting translations.
%To our knowledge, this is the first paper on the study of the nonlinear stability for the tumor model with vasculature.
\end{abstract}

\section{Introduction}
Over the past few decades, a variety of PDE models describing tumor growth in the form of free boundary problems have been proposed, developed and studied; see \cite{B,BC1,BC4,CLN,G1,G2,LFJC,WK1} and references therein. These models are based on mass conservation laws for cell densities and reaction-diffusion processes for nutrient concentrations within the tumor. Rigorous mathematical analysis and numerical simulation of these models has drawn considerable attention, and many interesting results have been established.
%in literatures \cite{BF,CF,C,C1,CE1,CE2,CF1,CF2,CF3,EM1,EM3,EM4,FH1,FH2,FH5,FL,FR1,FR3,HHH1,HHH2,HHH3,HZH1,HZH2,WC,WC3,ZEC1,ZW,ZWC} and references given there.
Mathematical analysis of such free boundary problems not only provides important insight into the growing mechanism of tumors, but may also aid in the understanding of experimental and clinical observations.

In this paper, we consider a free boundary tumor model with angiogenesis, consisting of proliferating cells only. This model was proposed by Friedman and Lam \cite{FL} as an extension to the model of Byrne and Chaplain \cite{BC1}.  Let $\Omega(t)$ denote the tumor region at time $t$, $\sigma$ and $p$ be the concentration of nutrients and pressure resulting from movement of cells, respectively. The equations are given by (see \cite{FL}):
 \begin{eqnarray}
   \sigma_t = \Delta\sigma-\sigma, \ \ \ \ \ \ \ \ x\in\Omega(t),\ t>0,\label{1.2} \\
   -\Delta p = \mu(\sigma-\widetilde{\sigma}), \ \ \ \ \ \ x\in\Omega(t),\ t>0, \label{1.4}
 \end{eqnarray}
where the positive parameter $\mu$ measures the aggressiveness of the tumor, and the positive constant $\widetilde{\sigma}$ is a threshold concentration for proliferation. The equation (\ref{1.4}) for the pressure is obtained by the conservation of mass (i.e., the cell proliferation rate $\mu(\sigma-\widetilde{\sigma})=\nabla\cdot\mathbf{V}$, here we assume a linear relationship between the cell proliferation rate and the nutrient concentration, $\mathbf{V}$ denotes the velocity field of the tumor cell movement) and the Darcy's law assumption (i.e., $\mathbf{V}=-\nabla p$, where the extracellular matrix is regarded as porous medium).

These equations are supplemented with the boundary conditions:
 \begin{eqnarray}
   \frac{\partial\sigma}{\partial \mathbf{n}}+\beta(\sigma-1) &=& 0, \ \ \ \ \ \ x\in\partial\Omega(t),\ t>0, \ \ \ \ \ (\widetilde{\sigma}<1),\label{1.3} \\
   p &=&  \gamma\kappa,  \ \ \ \ \ \ x\in\partial\Omega(t),\ t>0,\label{1.5}
 \end{eqnarray}
where in (\ref{1.3}) the impact from angiogenesis is incorporated (see \cite{FL}), $\mathbf{n}$ is the outward normal, (\ref{1.5}) describes the cell-to-cell adhesiveness (see \cite{B}) with $\gamma$ being the adhesiveness coefficient and $\kappa$ being the mean curvature.

If the velocity field is continuous up to the boundary,     the velocity of the free boundary is given by (see \cite{FL}):
\begin{equation}\label{1.6}
V_n=-\nabla p\cdot\mathbf{n}=-\frac{\partial p}{\partial \mathbf{n}},\ \ \ \ \ \  \ x\in\partial\Omega(t),\ t>0,
\end{equation}
where $V_\mathbf{n}$ is the velocity of the free boundary in the direction $\mathbf{n}$.

We finally impose the following initial conditions:
\begin{equation}\label{1.71}
  \sigma|_{t=0}=\sigma_0 \ \ \ \ \mbox{in}\ \ \Omega(0), \ \ \mbox{where\ } \Omega(0)\ \mbox{is\ given}.
\end{equation}

Angiogenesis is a process that tumor cells secret chemicals called tumor growth factors to stimulate the formation of new blood vessels penetrating into the tumor. In this model, nutrient enters the tumor only through these new blood vessels. The positive constant $\beta$ in (\ref{1.3}) reflects strength of the blood vessel system of the tumor:  $\beta=0$ means that the tumor does not have its own blood vessel system,  $\beta=\infty$ indicates that the tumor is all surrounded by the blood vessels  which reduces to the Dirichlet boundary condition $\sigma=1$.

The special case of this model $\beta=\infty$ corresponding the Dirichlet boundary condition $\sigma=1$ has drawn great attention of many researchers. Indeed, it was proved in \cite{FF} that a branch of symmetry-breaking stationary solutions bifurcates from the radially symmetric stationary solutions for each $\mu_n(R_S)\ (n\geq2)$ with free boundary
\begin{equation}\label{1.12}
r=R_S+\varepsilon Y_{n,0}(\theta)+O(\varepsilon^{2}).
\end{equation}
In \cite{BF} it was proved that if $\mu$ is sufficiently small, then under any small perturbations the stationary solution is asymptotically stable.
%In the sequel, Friedman and Hu \cite{FH1} determined the threshold value $\mu_\ast(R_S)$ for which the radial stationary solution
%is from linear stability to linear instability under small nonradial perturbations, and further found a critical number $\overline{R}$ ($\overline{R}\approx0.62207$ if $\overline{\sigma}=1$)
%such that $\mu_\ast(R_S)<\mu_2(R_S)$ if $R_S<\overline{R}$ and $\mu_\ast(R_S)=\mu_2(R_S)$ if $R_S>\overline{R}$. They also proved the spherical stationary solution is nonlinearly stable
%for $\mu(R_S)<\mu_\ast(R_S)$ in \cite{FH2}.
This work was improved by Friedman and Hu \cite{FH1,FH2}. They determined a threshold value $\mu^\ast$ for which the radially symmetric stationary solution changes from stability to instability under non-radial perturbations. Furthermore, linear stability analysis of the stationary solution at the bifurcation point $\mu=\mu^\ast$ was studied in \cite{FH5}. For the case that the nutrient consumption rate and the tumor cell proliferate rate are general functions, by using a center manifold analysis, Cui and Escher \cite{CE1,C3} found a positive critical value $\gamma^\ast$ (for the adhesiveness constant $\gamma$ in (\ref{1.5})) such that for $\gamma<\gamma^\ast$ the radially symmetric stationary solution is asymptotically stable with respect to small non-radial perturbations, while for $0<\gamma<\gamma^\ast$ this stationary solution is unstable.

Recently, Zhou and Wu \cite{ZW} considered the case that the nutrient concentration satisfies the Gibbs-Thomson relation on the boundary, and studied the existence and stability of the flat stationary solution. For the corresponding radially symmetric stationary solutions, it was proved \cite{W} that the stationary solution with larger radius is asymptotically stable and the other smaller one is unstable with respect to radial perturbations. In the sequel \cite{WZ}, they further refined the above result, and proved that under non-radial perturbations there exists a positive threshold value $\gamma^\ast$ such that the radially symmetric stationary solution with larger radius bifurcates from instability to stability.
 %is asymptotically stable for $\gamma>\gamma^\ast$ and unstable for $0<\gamma<\gamma^\ast$.
Moreover, considerable research works on asymptotic stability have been established for various tumor growth models, for instance, see \cite{C,CE2,EM1,WC3,ZH,ZH2,ZC,WC,FH4,WC1,WC2,CF1,C1,EM3,HHH2,HHH3,HZH3,ZX}.

%In the presence of inhibitor, the authors in \cite{WC} found a threshold value $\gamma^\ast$ for which the radially symmetric equilibrium changes from instability to stability with respect to non-radial perturbations and the asymptotic behavior of solutions was analyzed in \cite{C1,CF1}. If tumor growth is modeled to describe the growth of nonnecrotic tumors in different regimes of vascularization, the existence of the radially symmetric equilibrium (\hspace{-0.08em}\cite{EM4}) and its stability (\hspace{-0.08em}\cite{EM3}) were analyzed. Friedman, Cui and Chen \cite{CCF,C2,C1} investigated the tumor growth model describing the movement of various tumor cells (e.g., proliferating cells, quiescent cells, and dead cells) and established the asymptotic stability of the stationary solution. Numerical simulations on the stability of stationary solutions were carried out in \cite{HHH1,HHH2,HHH3}.

The radially symmetric version of the system (\ref{1.2})-(\ref{1.6}) was studied by Friedman and Lam \cite{FL}, in which they established the asymptotic behavior of the global solution and the existence of a unique radially symmetric stationary solution given by (see also \cite[Section 2]{HZH1})
\begin{equation}\label{1.7}
\begin{split}
    \sigma_S(r)=\frac{\beta}{\beta+R_SP_0(R_S)}\frac{R_S^{\frac{1}{2}}I_{1/2}(r)}{r^{\frac{1}{2}}I_{1/2}(R_S)}, \ \ \ 0<r<R_S,\ \ (P_n(R_S)\ \mbox{is\ defined\ in}\ (\ref{0.8})),
\end{split}
\end{equation}
\begin{equation}\label{1.70}
  p_S(r)=-\mu\sigma_S(r)+\frac{1}{6}\mu\widetilde{\sigma}r^{2}+C_1, \ \ \ 0<r<R_S,
\end{equation}
where
%\begin{equation}\label{1.8}
%\begin{split}
%    g(R_S)=\coth R_S-\frac{1}{R_S}=\frac{I_{3/2}(R_S)}{I_{1/2}(R_S)}=R_SP_0(R_S),\ \ (P_n(R_S)\ \mbox{is\ defined\ in}\ (\ref{0.8})),
%\end{split}
%\end{equation}
\begin{equation}\label{1.9}
    C_1=\frac{1}{R_S}+\frac{\mu\beta}{\beta+R_SP_0(R_S)}-\frac{1}{6}\mu\widetilde{\sigma}R_S^{2},
\end{equation}
and $R_S$ is the unique solution of
\begin{equation}\label{1.10}
    \frac{\beta P_0(R_S)}{\beta+R_SP_0(R_S)}=\frac{1}{3}\widetilde{\sigma}.
\end{equation}

Recently, it was established in \cite{HZH1} that $\mu_n(R_S)$ $(n\geq2\ \mbox{even})$ given by
\begin{equation}\label{1.11}
\begin{split}
    \mu_n(R_S)=\frac{3n(n-1)(n+2)}{2\widetilde{\sigma}R^4_S}\cdot\frac{\frac{n}{R_S}+\beta+R_SP_n(R_S)}{(n+\beta R_S)P_1(R_S)-(1+\beta R_S)P_n(R_S)}
    %&=\frac{\beta+R_SP_0(R_S)}{\beta R_SP_0(R_S)}\frac{n}{R^3_S}\left(\frac{n(n+1)}{2}-1\right)\frac{\frac{n}{R_S}+\beta+R_SP_n(R_S)}{(n+\beta R_S)P_1(R_S)-(1+\beta R_S)P_n(R_S)}
\end{split}
\end{equation}
 are bifurcation points of the symmetry-breaking stationary solution of the system (\ref{1.2})-(\ref{1.6}) with free boundary
\begin{equation*}
    r=R_S+\varepsilon Y_{n,0}(\theta)+o(\varepsilon),
\end{equation*}
where $Y_{n,0}$ is the spherical harmonic of order $(n,0)$. Moreover, $\mu_n(R_S)$ is monotonically increasing in $n$:
\begin{equation}\label{1.72}
  \mu_n(R_S)<\mu_{n+1}(R_S).
\end{equation}
In a recent paper \cite{HZH2}, we considered the linear stability of the stationary solution under non-radial perturbations and found a critical value $\mu^\ast=\mu^\ast(R_S)\ (\leq \mu_2(R_S))$ such that the stationary solution is stable for the linearized problem if $0<\mu<\mu^\ast(R_S)$, and it is unstable for the linearized problem if $\mu>\mu^\ast(R_S)$. With the linear stability of the radially symmetric stationary solution $(\sigma_S,p_S,R_S)$ already established, the present paper addresses the following question: {\em Is this stationary solution also stable for the original fully nonlinear free boundary problem for $0<\mu<\mu^\ast(R_S)$?}

Throughout the paper, by a rescaling if necessary we take $\gamma=1$. We assume the initial conditions are perturbed as follows:
\begin{equation}\label{1.13}
    \partial\Omega(0): \ r=R_S+\varepsilon\rho_0(\theta,\varphi), \ \ \ \ \sigma|_{t=0}=\sigma_S(r)+\varepsilon w_0(r,\theta,\varphi),
\end{equation}
where $\rho_0$ and $w_0$ are bounded functions.

\vskip 2mm
For the linearized problem, the translation of the origin resulting from mode 1 is easily seen. However, for the fully nonlinear problem, the perturbation of mode 1 is ``hidden" in the equation and it is not clear that mode 1 can be separated from other modes. Our grand challenge is to find a ``correct" translation of the origin to take care of perturbations from the mode 1 terms; this is characterized in the following (Theorem \ref{th2}) fixed point theorem.

The structure of this paper is as follows. In Section 2, we recall a few results from \cite{FH1,HZH2} which will be needed in the sequel. In Section 3, we transform the nonlinear problem into a new PDE system with a fixed boundary by Hanzawa transformation. In order to establish the stability result, it is crucial to find a new center of the sphere resulting from the perturbation of initial values; this is carried out in Section 4. This determination is crucial to our results. In Section 5, we derive the necessary estimates for all modes and complete the asymptotic stability in Section 6.

\section{Preliminaries}
In this section, we shall collect the properties of the function $P_n(\xi)$ that is introduced by Friedman and Hu \cite{FH1} and some results derived in the argument of the linear stability \cite{HZH2}, which are needed in our discussion.

The function $P_n(\xi)$ is defined by (see \cite{FH1})
\begin{equation}\label{0.8}
    P_n(\xi)=\frac{I_{n+3/2}(\xi)}{\xi I_{n+1/2}(\xi)}=2\sum_{m=1}^{\infty}\frac{1}{\xi^2+j_{n+1/2,m}^2}, \ \ \ n=0,1,2,3,\cdots
\end{equation}
where $j_{n,m}$ are the $m^{th}$ real positive zeros of Bessel functions $J_n(x)$, and $I_n(\xi)$ are the modified Bessel functions.

Recall \cite{FH1} that, for any $n\geq0$,
\begin{equation}\label{0.67}
  P_0(\xi)=\frac{1}{\xi}\coth\xi-\frac{1}{\xi^2},
\end{equation}
\begin{equation}\label{0.10}
    P_n(\xi)=\frac{1}{\xi^2P_{n+1}(\xi)+(2n+3)},
\end{equation}
%\begin{equation}\label{0.54}
%    \frac{d}{d\xi}P_n(\xi)=\frac{1}{\xi}-\frac{2n+3}{\xi}P_n(\xi)-\xi P_n^2(\xi),
%\end{equation}
%\begin{equation}\label{0.11}
%    P_n(0)=\frac{1}{2n+3},
%\end{equation}
%\begin{equation}\label{0.12}
%    \frac{d}{dr}\big\{rP_1(r)\big\}>0,
%\end{equation}
%\begin{equation}\label{0.13}
%    P_n(r)>P_{n+1}(r),
%\end{equation}
%\begin{equation}\label{0.14}
%    \frac{d}{dr}P_n(r)<0.
%\end{equation}
 %Also by \cite{FH1}, $P_n(\xi)$ can be rewritten as
%\begin{equation}\label{0.16}
%    P_n(\xi)=2\sum_{m=1}^{\infty}\frac{1}{\xi^2+j_{n+1/2,m}^2},
%\end{equation}
%where $j_{n,m}$ are the $m^{th}$ real positive zeros of $J_n(x)$ with $I_n(iz)=i^nJ_n(z)$.
\begin{equation}\label{6.14}
    \frac{d}{dr}\left(\frac{I_{n+1/2}(\sqrt{s+1}r)}{r^{1/2}}\right)=\frac{\sqrt{s+1}I_{n+3/2}(\sqrt{s+1}r)+\frac{n}{r}I_{n+1/2}(\sqrt{s+1}r)}{r^{1/2}}.
\end{equation}
As in \cite{HZH2}, we have
\begin{equation}\label{0.22}
    \lambda\triangleq\left(\frac{\partial^{2}\sigma_S}{\partial r^{2}}+\beta\frac{\partial\sigma_S}{\partial r}\right)\Big|_{r=R_S}=\frac{\beta P_0(R_S)}{\beta+R_SP_0(R_S)}[R_S^2P_1(R_S)+1+\beta R_S]
\end{equation}
and
\begin{equation}\label{0.56}
  \frac{\partial^2p_S}{\partial r^2}\Big|_{r=R}=-\mu\left(\frac{\beta}{\beta+R_SP_0(R_S)}-\widetilde{\sigma}\right).
\end{equation}
It is shown in \cite{HZH2} that the linear stability depends on the zeros of the following function:
\begin{equation}\label{1.34}
\begin{split}
    h_n(s)\triangleq h_n(s,\mu,R_S)=&\frac{1}{\mu}\frac{\beta+R_SP_0(R_S)}{\beta R_SP_0(R_S)}\left[s+\frac{n}{R_S^3}\left(\frac{n(n+1)}{2}-1\right)\right]-R_SP_1(R_S)\\
    &+\frac{R_S^2P_1(R_S)+1+\beta R_S}{(s+1)R_S+(\frac{n}{R_S}+\beta)/P_n(\sqrt{s+1}R_S)}.
\end{split}
\end{equation}
For $0<\mu<\mu^\ast=\mu^\ast(R_S)$, all zeros of $h_n(s)$ lie in $\mbox{Re}\ s<0$ ($n\neq1$); $s=0$ is a zero of $h_1(s)$, all other zeros lie in $\mbox{Re}\ s<0$ for $h_1(s)$. Furthermore,
\begin{equation}\label{0.036}
  \mu^\ast(R_S)<\mu_1(R_S),
\end{equation}
where
\begin{equation}\label{0.035}
    \mu_1(R_S)=\frac{\beta+R_SP_0(R_S)}{\beta R_S^3P_0(R_S)}\frac{R_S^2P_1(R_S)+1+\beta R_S}{R_SP_1^2(R_S)-\frac{1+\beta R_S}{2}P'_1(R_S)}.
\end{equation}
%and
%\begin{equation}\label{1.75}
%  \widetilde{\phi}_n(s)=(s+1)R+\Big(\frac{n}{R}+\beta\Big)/P_n(\sqrt{s+1}R), \ \ \ \ \phi_n(s)=P_n(\sqrt{s+1}R)\widetilde{\phi}_n(s).
%\end{equation}

As in the proof of \cite{CF}, we have the following local existence theorem:
\begin{thm}\label{th3}
If
\begin{equation}\label{0.57}
  (\sigma_0,\Gamma_0)\in C^{1+\gamma}(\overline{\Omega}(0))\times C^{4+\alpha}(\partial \Omega(0))\ \ \ \mbox{and}\ \ \frac{\partial\sigma_0}{\partial\mathbf{n}}+\beta(\sigma_0-1)=0 \ \ \ \mbox{on}\ \ \partial\Omega(0)
\end{equation}
for some $\alpha$, $\gamma$ $\in(0,1)$, then there exists a unique solution $(\sigma,p,\Gamma)$ of (\ref{1.2})-(\ref{1.6}) for $t\in[0,T]$ with some $T>0$, and
\begin{equation*}
  \sigma\in C^{1+\gamma,(1+\gamma)/2}\Big(\bigcup_{t\in[0,T]}\overline{\Omega}(t)\times\{t\}\Big)\cap C^{2+2\alpha/3,1+\alpha/3}\Big(\bigcup_{t\in[t_0,T]}\overline{\Omega}(t)\times\{t\}\Big), \ \ \ \mbox{for\ any}\ t_0>0,
\end{equation*}
\begin{equation*}
  p\in C^{2+\alpha,\alpha/3}\Big(\bigcup_{t\in[0,T]}\overline{\Omega}(t)\times\{t\}\Big), \ \ \ \Gamma\in C^{4+\alpha,1+\alpha/3}.
\end{equation*}
\end{thm}

%The following Laplace inverse transform formula is well known:
%\begin{lem}\label{5}
%Suppose that $F(s)$ is holomorphic in $\mbox{Re}\ s\geq a$ and
%\begin{equation*}
%    \limsup_{K\rightarrow\infty}\{|F(s)|; \mbox{Re}\ s\geq a,|s|\geq K\}=0.
%\end{equation*}
%Suppose also that there exists an $\eta>0$ and a constant $c_0$ such that
%\begin{equation*}
%    F(s)=c_0s^{-1}+O(|s|^{-1-\eta})\ \ \ \mbox{for}\ s=a+i\tau,\ \ \ \mbox{as}\ |\tau|\rightarrow\infty.
%\end{equation*}
%Then the integral
%\begin{equation*}
%    \psi(t)\equiv\frac{1}{2\pi i}\int^{a+i\infty}_{a-i\infty}F(s)e^{st}ds=\frac{e^{at}}{2\pi}\int^{\infty}_{-\infty}F(a+i\tau)e^{i\tau t}d\tau
%\end{equation*}
%defines the inverse Laplace transform of $F(s)$, i.e., $\widehat{g}(s)=F(s)$ for $\mbox{Re}\ s>a$, and furthermore,
%\begin{equation*}
%    |\psi(t)|\leq Ce^{at}.
%\end{equation*}
%\end{lem}

\section{The Nonlinearly Perturbed Problem}
The standard way to deal with free boundary problems is to transform it into a new PDE system with a fixed boundary. In this section, by Hanzawa transformation, we transform the free boundary problem (\ref{1.2})-(\ref{1.6}) into a nonlinearly perturbed problem in a fixed domain. For simplicity, we denote $R_S$ by $R$.

To begin with, let us assume the solution of the system (\ref{1.2})-(\ref{1.6}) is of the form (which will be verified in Section 6)
\begin{equation*}\label{3.1}
\begin{split}
    \partial\Omega(t):\ r&=R+\varepsilon\rho(\theta,\varphi,t),\\
    \sigma(r,\theta,\varphi,t)&=\sigma_S(r)+\varepsilon w(r,\theta,\varphi,t),\\
    p(r,\theta,\varphi,t)&=p_S(r)+\varepsilon q(r,\theta,\varphi,t),
\end{split}
\end{equation*}
and then the system (\ref{1.2})-(\ref{1.6}) can be written in terms of $(w,q,\rho)$ as follows:
%\begin{equation}\label{3.2}
%    w_t-\Delta w+w=0 \ \ \ \mbox{in} \ \Omega(t),\ t>0,
%\end{equation}
%\begin{equation}\label{3.3}
%    \frac{\partial w}{\partial \mathbf{n}}+\beta w=-\frac{1}{\varepsilon}\Big[\frac{\partial\sigma_S}{\partial\mathbf{n}}+\beta\sigma_S-\beta\Big]\ \ \ \mbox{on} \ \partial\Omega(t),\ t>0,
%\end{equation}
%\begin{equation}\label{3.4}
%    -\Delta q=\mu w \ \ \ \mbox{in} \ \Omega(t),\ t>0,
%\end{equation}
%\begin{equation}\label{3.5}
%    q=-\frac{1}{\varepsilon}\big[p_S-\kappa\big] \ \ \ \mbox{on} \ \partial\Omega(t),\ t>0,
%\end{equation}
%\begin{equation}\label{3.6}
%    \frac{d\rho}{dt}=-\Big(\frac{1}{\varepsilon}\frac{\partial p_S}{\partial\mathbf{n}}+\frac{\partial q}{\partial \mathbf{n}}\Big)\sqrt{1+\frac{|\varepsilon\nabla_\omega\rho|^2}{(R+\varepsilon\rho)^2}} \ \ \ \mbox{on} \ \partial\Omega(t),\ t>0,
%\end{equation}
\begin{eqnarray}
  w_t-\Delta w+w&=&0 \ \ \ \mbox{in} \ \Omega(t),\ t>0, \label{3.2}\\
  -\Delta q&=&\mu w \ \ \ \mbox{in} \ \Omega(t),\ t>0, \label{3.4}\\
  \frac{d\rho}{dt}&=&-\Big(\frac{1}{\varepsilon}\frac{\partial p_S}{\partial\mathbf{n}}+\frac{\partial q}{\partial \mathbf{n}}\Big)\sqrt{1+\frac{|\varepsilon\nabla_\omega\rho|^2}{(R+\varepsilon\rho)^2}} \ \ \ \mbox{on} \ \partial\Omega(t),\ t>0,\label{3.6}\\
  \frac{\partial w}{\partial \mathbf{n}}+\beta w&=&-\frac{1}{\varepsilon}\Big[\frac{\partial\sigma_S}{\partial\mathbf{n}}+\beta\sigma_S-\beta\Big]\ \ \ \mbox{on} \ \partial\Omega(t),\ t>0, \label{3.3}\\
  q&=&-\frac{1}{\varepsilon}\big[p_S-\kappa\big] \ \ \ \mbox{on} \ \partial\Omega(t),\ t>0, \label{3.5}
\end{eqnarray}
where $\nabla_\omega=\vec{e}_\theta\partial_\theta+\vec{e}_\varphi\frac{1}{\sin\theta}\partial_\varphi$.
%Note that after the initial values perturbed, say (\ref{1.13}), the system (\ref{3.2})-(\ref{3.6}) (or (\ref{1.2})-(\ref{1.6})) undergoes a translation of the origin resulting from the perturbation of mode 1. In order to prove asymptotic stability, we need to find the new center, which will be discussed in Section 4. We also note that the translation of the origin does not change the equation (\ref{3.2})-(\ref{3.6}).

Recall \cite{HZH1}, that
\begin{equation*}
  \mathbf{n}=\frac{1}{\sqrt{1+|\varepsilon\nabla_\omega\rho|^2/(R+\varepsilon\rho)^2}}
  \Big(\vec{e}_r-\frac{\varepsilon\rho_\theta}{R+\varepsilon\rho}\vec{e}_\theta-\frac{\varepsilon\rho_\varphi}{(R+\varepsilon\rho)\sin\theta}\vec{e}_\varphi\Big),
\end{equation*}
and
\begin{equation*}
 \nabla=\vec{e}_r\partial_r+\vec{e}_\theta\frac{1}{r}\partial_\theta+\vec{e}_\varphi\frac{1}{r\sin\theta}\partial_\varphi
 =\vec{e}_r\partial_r+\frac{1}{r}\nabla_\omega.
\end{equation*}
By the Taylor expansion, some of the right-hand sides terms of (\ref{3.6}) and (\ref{3.5}) can be written in the following way, respectively,
\begin{equation}\label{3.7}
\begin{split}
  \sqrt{1+\frac{|\varepsilon\nabla_\omega\rho|^2}{(R+\varepsilon\rho)^2}}\frac{\partial p_S(R+\varepsilon\rho)}{\partial\mathbf{n}}&=\sqrt{1+\frac{|\varepsilon\nabla_\omega\rho|^2}{(R+\varepsilon\rho)^2}}\cdot\nabla p_S|_{R+\varepsilon\rho}\mathbf{n}=\frac{\partial p_S(R+\varepsilon\rho)}{\partial r}\\
  &=\frac{\partial p_S(R)}{\partial r}+\frac{\partial^2 p_S(R)}{\partial r^2}\varepsilon\rho+\varepsilon^2 P_\varepsilon\\
  &=-\mu\left(\frac{\beta}{\beta+RP_0(R)}-\widetilde{\sigma}\right)\varepsilon\rho+\varepsilon^2 P_\varepsilon,\ \ (\mbox{by} \ (\ref{0.56})),
\end{split}
\end{equation}
%\begin{equation}\label{3.7}
%\begin{split}
%  &\sqrt{1+\frac{|\varepsilon\nabla_\omega\rho|^2}{(R+\varepsilon\rho)^2}}\Big(\frac{1}{\varepsilon}\frac{\partial p_S}{\partial\mathbf{n}}+\frac{\partial q}{\partial\mathbf{n}}\Big)\Big|_{R+\varepsilon\rho}\\
%  &\quad=\sqrt{1+\frac{|\varepsilon\nabla_\omega\rho|^2}{(R+\varepsilon\rho)^2}}\Big(\frac{1}{\varepsilon}\nabla p_S|_{R+\varepsilon\rho}\mathbf{n}+\nabla q|_{R+\varepsilon\rho}\mathbf{n}\Big)\\
%  &\quad=\frac{1}{\varepsilon}\frac{\partial p_S(R+\varepsilon\rho)}{\partial r}+\frac{\partial q}{\partial r}\Big|_{R+\varepsilon\rho}-\frac{\varepsilon}{(R+\varepsilon\rho)^2}\frac{\partial\rho}{\partial\theta}\frac{\partial q}{\partial\theta}
%  -\frac{\varepsilon}{(R+\varepsilon\rho)^2\sin^2\theta}\frac{\partial\rho}{\partial\varphi}\frac{\partial q}{\partial\varphi}\\
%  &\quad=\frac{\partial^2 p_S(R)}{\partial r^2}\rho+\frac{\partial q}{\partial r}\Big|_{R}+\varepsilon P_\varepsilon\\
%  &\quad=-\mu\left(\frac{\beta}{\beta+RP_0(R)}-\widetilde{\sigma}\right)\rho+\frac{\partial q}{\partial r}\Big|_{R}+\varepsilon P_\varepsilon,\ \ (\mbox{by} \ (\ref{0.56})),
%\end{split}
%\end{equation}
and
\begin{equation}\label{3.8}
  p_S(R+\varepsilon\rho)-\kappa=\frac{\varepsilon}{R^2}(\rho+\frac{1}{2}\Delta_\omega\rho)+\varepsilon^2K_\varepsilon,\ \ \ \mbox{(by\ \cite[Theorem 8.1]{FR2})},
\end{equation}
where
\begin{equation*}
    \Delta_\omega\rho=\frac{1}{\sin\theta}\frac{\partial}{\partial\theta}\left(\sin\theta\frac{\partial\rho}{\partial\theta}\right)
    +\frac{1}{\sin^2\theta}\frac{\partial^2\rho}{\partial\varphi^2}.
\end{equation*}
We now proceed to deal with (\ref{3.3}). In fact, (\ref{3.3}) is equal to
\begin{equation*}
\begin{split}
  &\frac{\partial w}{\partial r}\Big|_{R+\varepsilon\rho}-\frac{\varepsilon}{(R+\varepsilon\rho)^2}\frac{\partial\rho}{\partial\theta}\frac{\partial w}{\partial\theta}
  -\frac{\varepsilon}{(R+\varepsilon\rho)^2\sin^2\theta}\frac{\partial\rho}{\partial\varphi}\frac{\partial w}{\partial\varphi}+\beta w|_{R+\varepsilon\rho}\cdot\sqrt{1+\frac{|\varepsilon\nabla_\omega\rho|^2}{(R+\varepsilon\rho)^2}}\\
  &\quad=-\frac{1}{\varepsilon}\Big\{\frac{\partial\sigma_S(R+\varepsilon\rho)}{\partial r}
  +\beta\big[\sigma_S(R+\varepsilon\rho)-1\big]\cdot\sqrt{1+\frac{|\varepsilon\nabla_\omega\rho|^2}{(R+\varepsilon\rho)^2}}\Big\}
\end{split}
\end{equation*}
By the Taylor expansion, we have
\begin{equation}\label{3.9}
\begin{split}
  &\frac{\partial\sigma_S(R+\varepsilon\rho)}{\partial r}
  +\beta\big[\sigma_S(R+\varepsilon\rho)-1\big]\cdot\sqrt{1+\frac{|\varepsilon\nabla_\omega\rho|^2}{(R+\varepsilon\rho)^2}}\\
  &\quad=\frac{\partial\sigma_S(R+\varepsilon\rho)}{\partial r}
  +\beta\big[\sigma_S(R+\varepsilon\rho)-1\big]+\varepsilon^2\widetilde{S}_\varepsilon\\
  %&=-\frac{\partial^2 \sigma_S(R)}{\partial r^2}\rho-\beta\frac{\partial \sigma_S(R)}{\partial r}\rho+\varepsilon S_\varepsilon\\
  &\quad=\left(\frac{\partial^2 \sigma_S(R)}{\partial r^2}+\beta\frac{\partial \sigma_S(R)}{\partial r}\right)\varepsilon\rho+\varepsilon^2 S_\varepsilon\\
  &\quad=\varepsilon\lambda\rho+\varepsilon^2 S_\varepsilon,
\end{split}
\end{equation}

As in \cite{FH2}, by Hanzawa transformation which is defined by
\begin{equation*}
    r=r'+\chi(R-r')\varepsilon \rho(\theta,\varphi,t), \ \ \ t=t', \ \ \  \ \theta=\theta', \ \ \ \varphi=\varphi'
\end{equation*}
 with
\begin{equation*}
   \chi(z)\in C^{\infty}, \ \ \ \ \chi(z)=\left\{\begin{aligned}
&0,  \ \ \mbox{if} \ |z|\geq\frac{3}{4}\delta_0 \\
&1, \ \ \mbox{if}\ |z|<\frac{1}{4}\delta_0
\end{aligned}\right., \ \ \ \
 \left|\frac{d^{k}\chi}{dz^{k}}\right|\leq\frac{C}{\delta_0^{k}},\ \ \ \ (\delta_0\ \mbox{positive\ and\ small}),
\end{equation*}
 and the expansions (\ref{3.7})-(\ref{3.9}),
%of $\frac{\partial\sigma_S}{\partial\mathbf{n}}$, $\sigma_S$, $p_S$ and $\frac{\partial p_S}{\partial\mathbf{n}}$,
the system (\ref{3.2})-(\ref{3.5}) is transformed into the following system which is the nonlinearly perturbed problem of the system (\ref{1.2})-(\ref{1.6}) in a fixed domain,
\begin{equation}\label{3.10}
    w'_t-\Delta w'+w'=\varepsilon[-A_\varepsilon^1w'+A_\varepsilon w'] \ \ \ \mbox{in} \ B_R,\ t>0,
\end{equation}
\begin{equation}\label{3.12}
    \Delta' q'+\mu w'=-\varepsilon A_\varepsilon q' \ \ \ \mbox{in} \ B_R,\ t>0,
\end{equation}
\begin{equation}\label{3.14}
    \frac{d\rho'}{dt'}=\mu\left(\frac{\beta}{\beta+RP_0(R)}-\widetilde{\sigma}\right)\rho'-\frac{\partial q'}{\partial r'}+\varepsilon B^1_\varepsilon \ \ \ \mbox{on} \ \partial B_R,\ t>0,
\end{equation}
\begin{equation}\label{3.11}
    \frac{\partial w'}{\partial r'}+\beta w'=-\lambda\rho'+\varepsilon B^2_\varepsilon  \ \ \ \mbox{on} \ \partial B_R,\ t>0,
\end{equation}
\begin{equation}\label{3.13}
    q'=-\frac{1}{R^2}\Big(\rho'+\frac{1}{2}\Delta_\omega\rho'\Big)+\varepsilon B^3_\varepsilon \ \ \ \mbox{on} \ \partial B_R,\ t>0,
\end{equation}
where $B_R$ is the ball with radius $R$, $A_\varepsilon$ and $A_\varepsilon^1$ were given in \cite{FH2}, and
\begin{equation*}
\begin{split}
  B_\varepsilon^1=&-\frac{1}{\varepsilon^2}\frac{\partial p_S(R+\varepsilon\rho')}{\partial r'}-\frac{1}{\varepsilon}\mu\left(\frac{\beta}{\beta+RP_0(R)}-\widetilde{\sigma}\right)\rho'\\
  &+\frac{1}{(R+\varepsilon\rho')^2}\frac{\partial\rho'}{\partial\theta'}\frac{\partial q'}{\partial\theta'}
  +\frac{1}{(R+\varepsilon\rho')^2\sin^2\theta'}\frac{\partial\rho'}{\partial\varphi'}\frac{\partial q'}{\partial\varphi'},
\end{split}
\end{equation*}
\begin{equation*}
\begin{split}
  B^2_\varepsilon=&-\frac{1}{\varepsilon^2}\Big\{\frac{\partial\sigma_S(R+\varepsilon\rho')}{\partial r'}+\beta\big[\sigma_S(R+\varepsilon\rho')-1\big]\cdot\sqrt{1+\frac{|\varepsilon\nabla_\omega\rho'|^2}{(R+\varepsilon\rho')^2}}\Big\}\\
  &+\frac{1}{\varepsilon}\lambda\rho'+\frac{1}{\varepsilon}\Big\{-\beta w'\sqrt{1+\frac{|\varepsilon\nabla_\omega\rho'|^2}{(R+\varepsilon\rho')^2}}+\beta w'\Big\}\\
  &+\frac{1}{(R+\varepsilon\rho')^2}\frac{\partial\rho'}{\partial\theta'}\frac{\partial w'}{\partial\theta'}
  +\frac{1}{(R+\varepsilon\rho')^2\sin^2\theta'}\frac{\partial\rho'}{\partial\varphi'}\frac{\partial w'}{\partial\varphi'},
\end{split}
\end{equation*}
\begin{equation*}
  B_\varepsilon^3=-\frac{1}{\varepsilon^2}[p_S(R+\varepsilon\rho')-\kappa]+\frac{1}{\varepsilon R^2}\Big(\rho'+\frac{1}{2}\Delta_\omega\rho'\Big).
\end{equation*}

By the definition of $A_\varepsilon$ and $A_\varepsilon^1$ in \cite{FH2}, we obtain that $A_\varepsilon$ is a second order differential operator in $(r',\theta',\varphi')$, and $A_\varepsilon^1$ involves $\frac{\partial\rho}{\partial t}$ and a first order differential operator in $r'$. Furthermore, all terms of $A_\varepsilon^1$ and $A_\varepsilon$ do not involve any singularity, then it follows from Theorem \ref{th3} that, for $T>1$,
\begin{equation*}
  -A_\varepsilon^1w'+A_\varepsilon w'\in C^{2\alpha/3,\alpha/3}(B_R\times[0,T]),
  %\qquad \qquad A_\varepsilon q'\in C^{\alpha,\alpha/3}(B_R\times[0,T]).
\end{equation*}
\begin{equation*}
  A_\varepsilon q'\in C^{\alpha,\alpha/3}(B_R\times[0,T]).
\end{equation*}
Notice that the term $\frac{1}{\varepsilon}$ is cancelled out by the coefficient that accompanies it, so that $C^{2\alpha/3,\alpha/3}$ norm of $-A_\varepsilon^1w'+A_\varepsilon w'$ and $C^{\alpha,\alpha/3}$ norm of $A_\varepsilon q'$ are uniformly bounded in $\varepsilon$.

On the other hand, although $\sin^2\theta'$ appears in the denominator in the last term of $B_\varepsilon^1$, one can simply choose a different coordinate system to deal with this problem. As functions on $\Sigma=\{|x|=1\}$ (rather than as functions of $(\theta,\varphi)$), there are no singularities, and $B_\varepsilon^1\in C^{1+\alpha,\alpha/3}(\partial B_R\times[0,T])$. By (\ref{3.7}),
%\begin{equation*}
%\begin{split}
%  B_\varepsilon^1=&-\frac{1}{2}\frac{\partial^3p_S(R)}{\partial r'^3}\rho'^2+\frac{1}{(R+\varepsilon\rho')^2}\frac{\partial\rho'}{\partial\theta'}\frac{\partial q'}{\partial\theta'}\\
%  &+\frac{1}{(R+\varepsilon\rho')^2\sin^2\theta'}\frac{\partial\rho'}{\partial\varphi'}\frac{\partial q'}{\partial\varphi'}+O\Big(|\varepsilon|\Big\|\frac{\partial^4p_S(R)}{\partial r'^4}\Big\|\Big),
%\end{split}
%\end{equation*}
 $C^{1+\alpha,\alpha/3}$ norm of $B_\varepsilon^1$ is uniformly bounded in $\varepsilon$.
By (\ref{1.7}) and Theorem \ref{th3}, we obtain that $B_\varepsilon^2\in C^{1+2\alpha/3,1+\alpha/3}(\partial B_R\times[0,T])$. Moreover, we derive, by (\ref{3.9}),
%\begin{equation*}
%\begin{split}
%  B_\varepsilon^2=&-\frac{1}{2}\frac{\partial^3\sigma_S(R)}{\partial r'^3}\rho'^2-\frac{\beta}{2}\frac{\partial^2\sigma_S(R)}{\partial r'^2}\rho'^2+\frac{1}{2}\frac{\partial\sigma_S(R)}{\partial r'}\frac{|\nabla_\omega\rho'|^2}{(R+\varepsilon\rho')^2}\\
%  &+\frac{1}{(R+\varepsilon\rho')^2}\frac{\partial\rho'}{\partial\theta'}\frac{\partial w'}{\partial\theta'}+\frac{1}{(R+\varepsilon\rho')^2\sin^2\theta'}\frac{\partial\rho'}{\partial\varphi'}\frac{\partial w'}{\partial\varphi'}+O\Big(|\varepsilon|\Big\|\frac{\partial^4\sigma_S(R)}{\partial r'^4}\Big\|+|\varepsilon|\|\rho'\|_{C^2}\Big),
%\end{split}
%\end{equation*}
that $C^{1+2\alpha/3,1+\alpha/3}$ norm of $B_\varepsilon^2$ is uniformly bounded in $\varepsilon$.
Similarly, by (\ref{3.8}),
%\begin{equation*}
%\begin{split}
%  \kappa|_{R+\varepsilon\rho'}
%  %&=\frac{1}{R+\varepsilon\rho'}-\frac{\varepsilon}{2(R+\varepsilon\rho')^2}\Delta_\omega\rho'
%  %+O(\varepsilon^3\|\rho'\|_{C^2})\\
%  =\frac{1}{R}-\frac{\varepsilon}{R^2}\Big(\rho'+\frac{1}{2}\Delta_\omega\rho'\Big)+\frac{\varepsilon^2}{R^3}(\rho'^2+\rho'\Delta_\omega\rho')
%  +O(|\varepsilon|^3\|\rho'\|_{C^2}),
%\end{split}
%\end{equation*}
we find that $B_\varepsilon^3\in C^{2+\alpha,\alpha/3}(\partial B_R\times[0,T])$, and $C^{2+\alpha,\alpha/3}$ norm of $B_\varepsilon^3$ is uniformly bounded in $\varepsilon$.

 %By (\ref{1.7}), (\ref{1.70}) and Theorem \ref{th3}, $B_\varepsilon^1\in C^{1+\alpha,\alpha/3}(\partial B_R\times[0,T])$, $B_\varepsilon^2\in C^{2+2\alpha/3,1+\alpha/3}(\partial B_R\times[0,\infty))$, $B_\varepsilon^3\in C^{2+\alpha,\alpha/3}(\partial B_R\times[0,\infty))$. Furthermore, as can be seen by (\ref{3.9})-(\ref{3.7}), the above expressions $B_\varepsilon^1$, $B_\varepsilon^2$ and $B_\varepsilon^3$ are all uniformly bounded in $\varepsilon$.

For simplicity of notation, we shall denote functions $w'(r',\theta',\varphi',t')$, $q'(r',\theta',\varphi',t')$ and $\rho'(\theta',\varphi',t')$ again by $w(r,\theta,\varphi,t)$, $q(r,\theta,\varphi,t)$ and $\rho(\theta,\varphi,t)$, respectively, in the rest of this paper.

We need to estimate $(w,q,\rho)$ of the system (\ref{3.10})-(\ref{3.13}) to prove asymptotic stability. However the system (\ref{3.10})-(\ref{3.13}) is nonlinear, the method we shall use is to study the inhomogeneous linear system instead of the nonlinearly perturbed system (\ref{3.10})-(\ref{3.13}), namely, we consider the system (\ref{3.10})-(\ref{3.13}) where the $\varepsilon$ terms of any order are replaced by given functions,
\begin{equation}\label{3.15}
    w_t-\Delta w+w=\varepsilon f^1(r,\theta,\varphi,t) \ \ \ \mbox{in} \ B_R,\ t>0,
\end{equation}
\begin{equation}\label{3.17}
    \Delta q+\mu w=-\varepsilon f^2(r,\theta,\varphi,t) \ \ \ \mbox{in} \ B_R,\ t>0,
\end{equation}
\begin{equation}\label{3.19}
    \frac{d\rho}{dt}=\mu\left(\frac{\beta}{\beta+RP_0(R)}-\widetilde{\sigma}\right)\rho-\frac{\partial q}{\partial r}+\varepsilon b^1(\theta,\varphi,t) \ \ \ \mbox{on} \ \partial B_R,\ t>0,
\end{equation}
\begin{equation}\label{3.16}
    \frac{\partial w}{\partial r}+\beta w=-\lambda\rho+\varepsilon b^2(\theta,\varphi,t)  \ \ \ \mbox{on} \ \partial B_R,\ t>0,
\end{equation}
\begin{equation}\label{3.18}
    q=-\frac{1}{R^2}\Big(\rho+\frac{1}{2}\Delta_\omega\rho\Big)+\varepsilon b^3(\theta,\varphi,t) \ \ \ \mbox{on} \ \partial B_R,\ t>0,
\end{equation}
with initial conditions
\begin{equation}\label{3.20}
  w|_{t=0}=w_0, \ \ \ \quad \ \ \ \rho|_{t=0}=\rho_0.
\end{equation}
%\begin{equation}\label{3.21}
%  \rho|_{t=0}=\rho_0.
%\end{equation}
Then for given functions $f^i$, $b^j$ satisfying additional assumptions, we solve the inhomogeneous linear system (\ref{3.15})-(\ref{3.18}), and derive the estimate of $(w,q,\rho)$. After that we want to define the new functions $\widetilde{f}^i$, $\widetilde{b}^j$ by
\begin{equation*}
  \widetilde{f}^1=-A_\varepsilon^1w+A_\varepsilon w, \ \ \ \ \widetilde{f}^2=-A_\varepsilon^1w,\ \ \ \
  \widetilde{b}^1=B_\varepsilon^1, \ \ \ \widetilde{b}^2=B_\varepsilon^2, \ \ \ \widetilde{b}^3=B_\varepsilon^3.
\end{equation*}
We shall show that the mapping $S:(f^i,b^j)\rightarrow(\widetilde{f}^i,\widetilde{b}^j)$ admits a fixed point.
 %for $t\in[0,T]$ with some $T>0$ and then establish the asymptotic estimate for the system (\ref{1.2})-(\ref{1.6}) on the interval $[0,T]$, and then proceed to establish the same estimate for all $[0,\infty)$.

Assume that the functions $f^i$, $b^j$ satisfy
\begin{equation}\label{4.1}
  \sqrt{|\varepsilon|}\left(\int_0^\infty e^{2\delta_1t}\|f^1(\cdot,t)\|^2_{L^2(B_R)}dt\right)^{1/2}\leq1,
\end{equation}
\begin{equation}\label{4.2}
  \sqrt{|\varepsilon|}\left(\int_0^\infty e^{2\delta_1t}\|f^2(\cdot,t)\|^2_{L^2(B_R)}dt\right)^{1/2}\leq1,
\end{equation}
\begin{equation}\label{4.3}
  \sqrt{|\varepsilon|}\left(\int_0^\infty e^{2\delta_1t}\|b^1(\cdot,t)\|^2_{H^{1/2}(\partial B_R)}dt\right)^{1/2}\leq1,
\end{equation}
\begin{equation}\label{4.4}
  \sqrt{|\varepsilon|}\left(\int_0^\infty e^{2\delta_1t}\|b^2(\cdot,t)\|^2_{H^{1/2}(\partial B_R)}dt\right)^{1/2}\leq1,
\end{equation}
\begin{equation}\label{4.5}
  \sqrt{|\varepsilon|}\left(\int_0^\infty e^{2\delta_1t}\|b^3(\cdot,t)\|^2_{H^{3/2}(\partial B_R)}dt\right)^{1/2}\leq1,
\end{equation}
where $\delta_1$ is positive and sufficiently small, and that, for some $\alpha\in(0,1)$,
\begin{equation}\label{4.6}
  \sqrt{|\varepsilon|}\|f^1\|_{C^{2\alpha/3,\alpha/3}(B_R\times[0,\infty))}\leq1,
\end{equation}
\begin{equation}\label{4.7}
  \sqrt{|\varepsilon|}\|f^2\|_{C^{\alpha,\alpha/3}(B_R\times[0,\infty))}\leq1,
\end{equation}
\begin{equation}\label{4.8}
  \sqrt{|\varepsilon|}\|b^1\|_{C^{1+\alpha,\alpha/3}(\partial B_R\times[0,\infty))}\leq1,
\end{equation}
\begin{equation}\label{4.9}
  \sqrt{|\varepsilon|}\|b^2\|_{C^{1+2\alpha/3,1+\alpha/3}(\partial B_R\times[0,\infty))}\leq1,
\end{equation}
\begin{equation}\label{4.10}
  \sqrt{|\varepsilon|}\|b^3\|_{C^{2+\alpha,\alpha/3}(\partial B_R\times[0,\infty))}\leq1.
\end{equation}
For simplicity, we also assume that
\begin{equation}\label{4.12}
  \rho_0\in C^{4+\alpha}(\overline{B}_R), \ \ \ \ \ w_0\in C^{2+2\alpha/3}(\overline{B}_R),
\end{equation}
and the compatibility condition of order 2 for $w$ is satisfied.

We formally expand all the functions $f^i$, $b^j$ in terms of spherical harmonics
\begin{equation*}
  f^i(r,\theta,\varphi,t)=\sum^\infty_{n=0}\sum^n_{m=-n}f^i_{n,m}(r,t)Y_{n,m}(\theta,\varphi), \ \ \ \ i=1,2,
\end{equation*}
\begin{equation*}
  b^j(\theta,\varphi,t)=\sum^\infty_{n=0}\sum^n_{m=-n}b^j_{n,m}(t)Y_{n,m}(\theta,\varphi), \ \ \ \ j=1,2,3.
\end{equation*}
Recall \cite{FH2} that (\ref{4.1})-(\ref{4.5}) imply
\begin{equation}\label{4.19}
  |\varepsilon|\int_0^\infty e^{2\delta_1t}\|f_{n,m}^1(\cdot,t)\|^2_{L^2(B_R)}dt=F_{n,m}^1, \ \ \ \ \ \sum_{n,m}F_{n,m}^1\leq1,
\end{equation}
\begin{equation}\label{4.20}
  |\varepsilon|\int_0^\infty e^{2\delta_1t}\|f_{n,m}^2(\cdot,t)\|^2_{L^2(B_R)}dt=F_{n,m}^2, \ \ \ \ \ \sum_{n,m}F_{n,m}^2\leq1,
\end{equation}
\begin{equation}\label{4.21}
  |\varepsilon|(n+1)\int_0^\infty e^{2\delta_1t}|b_{n,m}^1(t)|^2dt=B_{n,m}^1, \ \ \ \ \ \sum_{n,m}B_{n,m}^1\leq C,
\end{equation}
\begin{equation}\label{4.22}
  |\varepsilon|(n+1)\int_0^\infty e^{2\delta_1t}|b_{n,m}^2(t)|^2dt=B_{n,m}^2, \ \ \ \ \ \sum_{n,m}B_{n,m}^2\leq C,
\end{equation}
\begin{equation}\label{4.23}
  |\varepsilon|(n+1)^3\int_0^\infty e^{2\delta_1t}|b_{n,m}^3(t)|^2dt=B_{n,m}^3, \ \ \ \ \ \sum_{n,m}B_{n,m}^3\leq C.
\end{equation}

Now we proceed to solve the inhomogeneous linear system (\ref{3.15})-(\ref{3.18}). We look for a solution of the form:
\begin{eqnarray*}
  w(r,\theta,\varphi,t) &=& \sum^\infty_{n=0}\sum^n_{m=-n}w_{n,m}(r,t)Y_{n,m}(\theta,\varphi), \\
   q(r,\theta,\varphi,t) &=& \sum^\infty_{n=0}\sum^n_{m=-n}q_{n,m}(r,t)Y_{n,m}(\theta,\varphi), \\
\rho(\theta,\varphi,t) &=& \sum^\infty_{n=0}\sum^n_{m=-n}\rho_{n,m}(t)Y_{n,m}(\theta,\varphi).
\end{eqnarray*}
Then, by the relation
\begin{equation*}\label{1.20}
    \Delta_\omega Y_{n,m}(\theta,\varphi)+n(n+1)Y_{n,m}(\theta,\varphi)=0\ \ \ \mbox{and}\ \ \ \Delta=\frac{1}{r^2}\frac{\partial}{\partial r}\Big(r^2\frac{\partial}{\partial r}\Big)+\frac{1}{r^2}\Delta_\omega,
\end{equation*}
$w_{n,m}(r,t)$, $q_{n,m}(r,t)$ and $\rho_{n,m}(t)$ satisfy the following inhomogeneous linear system
\begin{equation}\label{4.13}
    \partial_tw_{n,m}(r,t)-\Delta w_{n,m}(r,t)+\left(\frac{n(n+1)}{r^2}+1\right)w_{n,m}(r,t)=\varepsilon f^1_{n,m}(r,t) \ \ \ \mbox{in}\ B_R\times\{ t>0\},
\end{equation}
\begin{equation}\label{4.15}
    -\Delta q_{n,m}(r,t)+\frac{n(n+1)}{r^2}q_{n,m}(r,t)=\mu w_{n,m}(r,t)+\varepsilon f^2_{n,m}(r,t) \ \ \ \mbox{in} \ B_R\times\{ t>0\},
\end{equation}
\begin{equation}\label{4.17}
    \frac{d\rho_{n,m}(t)}{dt}=\mu\left(\frac{\beta}{\beta+RP_0(R)}-\widetilde{\sigma}\right)\rho_{n,m}(t)-\frac{\partial q_{n,m}(R,t)}{\partial r}+\varepsilon b^1_{n,m}(t), \ \ \ t>0,
\end{equation}
\begin{equation}\label{4.14}
    \frac{\partial w_{n,m}(R,t)}{\partial r}+\beta w_{n,m}(R,t)=-\lambda\rho_{n,m}(t)+\varepsilon b^2_{n,m}(t),  \ \ \  \ \ t>0,
\end{equation}
\begin{equation}\label{4.16}
    q_{n,m}(R,t)=-\frac{1}{R^2}\left(1-\frac{n(n+1)}{2}\right)\rho_{n,m}(t)+\varepsilon b^3_{n,m}(t), \ \ \  \ \ t>0,
\end{equation}
\begin{equation}\label{4.18}
    \rho_{n,m}|_{t=0}=\rho_{0,n,m}, \ \ \ \ \ w_{n,m}|_{t=0}=w_{0,n,m}(r) \ \ \ \ \mbox{in}\ \ B_R.
\end{equation}

Introduce the Laplace transform,
%\begin{equation*}
%\begin{split}
  $$\widehat{w}_{n,m}(r,s)=\int^{\infty}_0e^{-st}w_{n,m}(r,t)dt, \ \ \ \widehat{q}_{n,m}(r,s) = \int^{\infty}_0e^{-st}q_{n,m}(r,t)dt,$$
   $$\widehat{\rho}_{n,m}(s) = \int^{\infty}_0e^{-st}\rho_{n,m}(t)dt.$$
%\end{split}
%\end{equation*}
Taking formally the Laplace transform of (\ref{4.13})-(\ref{4.16}),
 we get
%\begin{eqnarray}
%% \nonumber to remove numbering (before each equation)
%  -\Delta \widehat{w}(r,s)+(\frac{n(n+1)}{r^2}+(s+1))\widehat{w}(r,s) &=& w_{0,n,m}(r)+\varepsilon \widehat{f}^1_{n,m}(r,s) \ \ \mbox{in} \ B_R,\label{4.24} \\
%  \frac{\partial \widehat{w}(R,s)}{\partial r}+\beta \widehat{w}(R,s) &=& -\lambda\widehat{\rho}+\varepsilon \widehat{b}^1_{n,m}(s),\label{4.25} \\
%  -\Delta \widehat{q}(r,s)+\frac{n(n+1)}{r^2}\widehat{q}(r,s) &=& \mu \widehat{w}(r,s)+\varepsilon \widehat{f}^2_{n,m}(r,s) \ \ \ \ \  \mbox{in} \ B_R,\label{4.26} \\
%  \widehat{q}(R,s) &=& \frac{1}{R^2}\left(\frac{n(n+1)}{2}-1\right)\widehat{\rho}(s)+\varepsilon \widehat{b}^2_{n,m}(s),\label{4.27} \\
%  s\widehat{\rho}(s)-\rho_{0,n,m} &=& \mu\left(\frac{\beta}{\beta+RP_0(R)}-\widetilde{\sigma}\right)\widehat{\rho}(s)-\frac{\partial \widehat{q}}{\partial r}\Big|_{r=R}+\varepsilon \widehat{b}^3_{n,m}(s).\label{4.28}
%\end{eqnarray}
\begin{equation}\label{4.24}
    -\Delta \widehat{w}_{n,m}(r,s)+\Big(\frac{n(n+1)}{r^2}+s+1\Big)\widehat{w}_{n,m}(r,s) = w_{0,n,m}(r)+\varepsilon \widehat{f}^1_{n,m}(r,s) \ \ \mbox{in} \ B_R,
\end{equation}
\begin{equation}\label{4.26}
    -\Delta \widehat{q}_{n,m}(r,s)+\frac{n(n+1)}{r^2}\widehat{q}_{n,m}(r,s) = \mu \widehat{w}_{n,m}(r,s)+\varepsilon \widehat{f}^2_{n,m}(r,s) \ \ \ \ \  \mbox{in} \ B_R,
\end{equation}
\begin{equation}\label{4.28}
    s\widehat{\rho}_{n,m}(s)-\rho_{0,n,m}=\mu\left(\frac{\beta}{\beta+RP_0(R)}-\widetilde{\sigma}\right)\widehat{\rho}_{n,m}(s)
    -\frac{\partial \widehat{q}_{n,m}}{\partial r}\Big|_{r=R}+\varepsilon \widehat{b}^1_{n,m}(s).
\end{equation}
\begin{equation}\label{4.25}
   \frac{\partial \widehat{w}_{n,m}(R,s)}{\partial r}+\beta \widehat{w}_{n,m}(R,s) = -\lambda\widehat{\rho}_{n,m}(s)+\varepsilon \widehat{b}^2_{n,m}(s),
\end{equation}
\begin{equation}\label{4.27}
   \widehat{q}_{n,m}(R,s) = \frac{1}{R^2}\left(\frac{n(n+1)}{2}-1\right)\widehat{\rho}_{n,m}(s)+\varepsilon \widehat{b}^3_{n,m}(s).
\end{equation}

As in \cite{HZH2}, using (\ref{6.14}), we can solve (\ref{4.24}) and (\ref{4.25}) in the form
\begin{equation}\label{4.29}
\begin{split}
    \widehat{w}_{n,m}(r,s)=&\frac{(-\lambda\widehat{\rho}_{n,m}(s)+\varepsilon \widehat{b}^2_{n,m}(s))R^{\frac{1}{2}}}{\sqrt{s+1}I_{n+3/2}(\sqrt{s+1}R)+(\frac{n}{R}+\beta) I_{n+1/2}(\sqrt{s+1}R)}\frac{I_{n+1/2}(\sqrt{s+1}r)}{r^{\frac{1}{2}}}\\
    &+\xi_{1,n,m}(r,s)+\varepsilon\xi_{2,n,m}(r,s),
\end{split}
\end{equation}
where $\xi_{1,n,m}(r,s)$ is the solution of
\begin{equation}\label{4.30}
\begin{split}
  -\Delta\xi_{1,n,m}(r,s)+\left(\frac{n(n+1)}{r^2}+s+1\right)\xi_{1,n,m}(r,s) &= w_{0,n,m}(r) \ \ \ \ \mbox{in}\ B_R, \\
  \left(\frac{\partial \xi_{1,n,m}(r,s)}{\partial r}+\beta \xi_{1,n,m}(r,s)\right)\Big|_{r=R} &= 0,
\end{split}
\end{equation}
and $\xi_{2,n,m}(r,s)$ is the solution of
\begin{equation}\label{4.300}
\begin{split}
  -\Delta\xi_{2,n,m}(r,s)+\left(\frac{n(n+1)}{r^2}+s+1\right)\xi_{2,n,m}(r,s) &=\widehat{f}^1_{n,m}(r,s)  \ \ \ \ \mbox{in}\ B_R, \\
  \left(\frac{\partial \xi_{2,n,m}(r,s)}{\partial r}+\beta \xi_{2,n,m}(r,s)\right)\Big|_{r=R} &= 0.
\end{split}
\end{equation}

Let
\begin{equation}\label{4.31}
    \phi=\widehat{q}_{n,m}+\frac{\mu}{s+1}\widehat{w}_{n,m},
\end{equation}
then $\phi$ satisfies
\begin{equation}\label{4.32}
    -\Delta\phi+\frac{n(n+1)}{r^2}\phi=\frac{\mu}{s+1}(w_{0,n,m}+\varepsilon \widehat{f}^1_{n,m})+\varepsilon \widehat{f}^2_{n,m} \ \ \ \ \mbox{in}\ B_R,
\end{equation}
and by (\ref{4.28}), (\ref{4.27}) and (\ref{4.25}), $\phi$ satisfies the following boundary condition
\begin{equation}\label{4.33}
\begin{split}
    &\Big(\frac{\partial\phi}{\partial r}+\beta\phi\Big)\Big|_{r=R}\\
    &=\mu\left(\frac{\beta}{\beta+RP_0(R)}-\widetilde{\sigma}\right)\widehat{\rho}_{n,m}(s)-(s\widehat{\rho}_{n,m}(s)-\rho_{0,n,m})+\varepsilon \widehat{b}^1_{n,m}(s)\\
    &\quad+\frac{\beta}{R^2}\left(\frac{n(n+1)}{2}-1\right)\widehat{\rho}_{n,m}(s)+\beta\varepsilon \widehat{b}^3_{n,m}(s)-\frac{\mu\lambda}{s+1}\widehat{\rho}_{n,m}(s)+\frac{\mu}{s+1}\varepsilon \widehat{b}^2_{n,m}(s).
\end{split}
\end{equation}
The solution of the problem (\ref{4.32})-(\ref{4.33}) is given by
\begin{equation}\label{4.34}
\begin{split}
    \phi(r,s)=&\frac{r^n}{nR^{n-1}+\beta R^n}\left[\mu\left(\frac{\beta}{\beta+RP_0(R)}-\widetilde{\sigma}\right)\widehat{\rho}_{n,m}(s)-(s\widehat{\rho}_{n,m}(s)-\rho_{0,n,m})\right.\\
    &\left.+\frac{\beta}{R^2}\left(\frac{n(n+1)}{2}-1\right)\widehat{\rho}_{n,m}(s)-\frac{\mu\lambda}{s+1}\widehat{\rho}_{n,m}(s)+\varepsilon \widehat{b}^1_{n,m}(s)\right.\\
    &\left.+\frac{\mu}{s+1}\varepsilon \widehat{b}^2_{n,m}(s)+\beta\varepsilon \widehat{b}^3_{n,m}(s)\right]+\phi_{1,n,m}(r,s)+\varepsilon\phi_{2,n,m}(r,s),
\end{split}
\end{equation}
where $\phi_{1,n,m}(r,s)$ is the solution of
\begin{equation}\label{4.35}
\begin{split}
  -\Delta\phi_{1,n,m}+\frac{n(n+1)}{r^2}\phi_{1,n,m} &= \frac{\mu}{s+1}w_{0,n,m}\ \ \ \ \mbox{in}\ B_R, \\
  \left(\frac{\partial\phi_{1,n,m}}{\partial r}+\beta\phi_{1,n,m}\right)\Big|_{r=R} &= 0,
\end{split}
\end{equation}
and $\phi_{2,n,m}(r,s)$ is the solution of
\begin{equation}\label{4.36}
\begin{split}
  -\Delta\phi_{2,n,m}+\frac{n(n+1)}{r^2}\phi_{2,n,m} &= \frac{\mu}{s+1}\widehat{f}^1_{n,m}+\widehat{f}^2_{n,m}\ \ \ \ \mbox{in}\ B_R, \\
  \left(\frac{\partial\phi_{2,n,m}}{\partial r}+\beta\phi_{2,n,m}\right)\Big|_{r=R} &= 0.
\end{split}
\end{equation}

By (\ref{4.31}), (\ref{4.34}), (\ref{4.29}) and (\ref{6.14}), we have
\begin{equation*}\label{4.37}
\begin{split}
    \frac{\partial\widehat{q}_{n,m}}{\partial r}\Big|_{r=R}=&\frac{\partial\phi}{\partial r}\Big|_{r=R}-\frac{\mu}{s+1}\frac{\partial\widehat{w}_{n,m}}{\partial r}\Big|_{r=R}\\
    =&\frac{n}{n+\beta R}\Big\{\mu\Big(\frac{\beta}{\beta+RP_0(R)}-\widetilde{\sigma}\Big)\widehat{\rho}_{n,m}(s)-(s\widehat{\rho}_{n,m}(s)-\rho_{0,n,m})\\
    &+\frac{\beta}{R^2}\Big(\frac{n(n+1)}{2}-1\Big)\widehat{\rho}_{n,m}(s)-\frac{\mu\lambda}{s+1}\widehat{\rho}_{n,m}(s)\\
    &+\varepsilon \widehat{b}^1_{n,m}(s)+\frac{\mu}{s+1}\varepsilon \widehat{b}^2_{n,m}(s)+\beta\varepsilon \widehat{b}^3_{n,m}(s)\Big\}+\frac{\partial\phi_{1,n,m}}{\partial r}\Big|_{r=R}+\varepsilon\frac{\partial\phi_{2,n,m}}{\partial r}\Big|_{r=R}\\
    &-\frac{\mu}{s+1}\Big\{\frac{(-\lambda\widehat{\rho}_{n,m}(s)+\varepsilon\widehat{b}^2_{n,m})[\sqrt{s+1}I_{n+3/2}(\sqrt{s+1}R)+\frac{n}{R} I_{n+1/2}(\sqrt{s+1}R)]}{\sqrt{s+1}I_{n+3/2}(\sqrt{s+1}R)+(\frac{n}{R}+\beta) I_{n+1/2}(\sqrt{s+1}R)}\\
    &+\frac{\partial\xi_{1,n,m}}{\partial r}\Big|_{r=R}+\varepsilon\frac{\partial\xi_{2,n,m}}{\partial r}\Big|_{r=R}\Big\}.
\end{split}
\end{equation*}
Inserting the above expression into (\ref{4.28}), we obtain
\begin{equation*}
\begin{split}
   &-\frac{\beta R}{n+\beta R}\mu\Big(\frac{\beta}{\beta+RP_0(R)}-\widetilde{\sigma}\Big)\widehat{\rho}_{n,m}(s)+\frac{\beta R}{n+\beta R}(s\widehat{\rho}_{n,m}(s)-\rho_{0,n,m})\\
   &+\frac{n}{n+\beta R}\frac{\beta}{R^2}\Big(\frac{n(n+1)}{2}-1\Big)\widehat{\rho}_{n,m}(s)-\frac{n}{n+\beta R}\frac{\mu\lambda}{s+1}\widehat{\rho}_{n,m}(s)\\
   &+\varepsilon\frac{n}{n+\beta R}\big[\widehat{b}^1_{n,m}+\frac{\mu}{s+1}\widehat{b}^2_{n,m}+\beta\widehat{b}^3_{n,m}\big]+\frac{\partial\phi_{1,n,m}}{\partial r}\Big|_{r=R}+\varepsilon\frac{\partial\phi_{2,n,m}}{\partial r}\Big|_{r=R}\\
   &-\frac{\mu}{s+1}(-\lambda\widehat{\rho}_{n,m}(s)+\varepsilon\widehat{b}^2_{n,m})\frac{\sqrt{s+1}I_{n+3/2}(\sqrt{s+1}R)+\frac{n}{R} I_{n+1/2}(\sqrt{s+1}R)}{\sqrt{s+1}I_{n+3/2}(\sqrt{s+1}R)+(\frac{n}{R}+\beta) I_{n+1/2}(\sqrt{s+1}R)}\\
   &-\frac{\mu}{s+1}\Big(\frac{\partial\xi_{1,n,m}}{\partial r}\Big|_{r=R}+\varepsilon\frac{\partial\xi_{2,n,m}}{\partial r}\Big|_{r=R}\Big)-\varepsilon\widehat{b}^1_{n,m}=0.
\end{split}
\end{equation*}
After a direct computation, this equality is equivalent to
\begin{equation*}
\begin{split}
   &\Big\{\frac{\beta R}{n+\beta R}s+\frac{n}{n+\beta R}\frac{\beta}{R^2}\Big(\frac{n(n+1)}{2}-1\Big)-\frac{\beta R}{n+\beta R}\mu\Big(\frac{\beta}{\beta+RP_0(R)}-\widetilde{\sigma}\Big)\\
    &\quad-\frac{n}{n+\beta R}\frac{\mu\lambda}{s+1}+\frac{\mu\lambda}{s+1}\cdot\frac{\sqrt{s+1}I_{n+3/2}(\sqrt{s+1}R)+\frac{n}{R} I_{n+1/2}(\sqrt{s+1}R)}{\sqrt{s+1}I_{n+3/2}(\sqrt{s+1}R)+(\frac{n}{R}+\beta) I_{n+1/2}(\sqrt{s+1}R)}\Big\}\widehat{\rho}_{n,m}(s)\\
    &=\frac{\beta R}{n+\beta R}\rho_{0,n,m}-\varepsilon\frac{n}{n+\beta R}\big[\widehat{b}^1_{n,m}+\frac{\mu}{s+1}\widehat{b}^2_{n,m}+\beta\widehat{b}^3_{n,m}\big]\\
    &\quad+\varepsilon\frac{\mu}{s+1}\widehat{b}^2_{n,m}\frac{\sqrt{s+1}I_{n+3/2}(\sqrt{s+1}R)+\frac{n}{R} I_{n+1/2}(\sqrt{s+1}R)}{\sqrt{s+1}I_{n+3/2}(\sqrt{s+1}R)+(\frac{n}{R}+\beta) I_{n+1/2}(\sqrt{s+1}R)}\\
    &\quad-\Big(\frac{\partial\phi_{1,n,m}}{\partial r}\Big|_{r=R}+\varepsilon\frac{\partial\phi_{2,n,m}}{\partial r}\Big|_{r=R}\Big)+\frac{\mu}{s+1}\Big(\frac{\partial\xi_{1,n,m}}{\partial r}\Big|_{r=R}+\varepsilon\frac{\partial\xi_{2,n,m}}{\partial r}\Big|_{r=R}\Big)+\varepsilon\widehat{b}^1_{n,m}.\\
\end{split}
\end{equation*}
By \cite[page 2482]{HZH2} and (\ref{1.34}), the brace in left-hand side is equal to
$$\frac{\beta R}{n+\beta R}\frac{\beta RP_0(R)}{\beta+RP_0(R)}\mu h_n(s,\mu,R),$$
and the right-hand side is equal to
\begin{equation*}\label{0.53}
\begin{split}
    %&\frac{\beta R}{n+\beta R}\mu\frac{\beta RP_0(R)}{\beta+RP_0(R)}h_n(s,\mu,R)\widehat{\rho}(s)=\\
    %&\Big\{\frac{\beta R}{n+\beta R}s+\frac{n}{n+\beta R}\frac{\beta}{R^2}\Big(\frac{n(n+1)}{2}-1\Big)-\frac{\beta R}{n+\beta R}\mu\Big(\frac{\beta}{\beta+RP_0(R)}-\widetilde{\sigma}\Big)\\
%    &\quad-\frac{n}{n+\beta R}\frac{\mu\lambda}{s+1}+\frac{\mu\lambda}{s+1}\cdot\frac{\sqrt{s+1}I_{n+3/2}(\sqrt{s+1}R)+\frac{n}{R} I_{n+1/2}(\sqrt{s+1}R)}{\sqrt{s+1}I_{n+3/2}(\sqrt{s+1}R)+(\frac{n}{R}+\beta) I_{n+1/2}(\sqrt{s+1}R)}\Big\}\widehat{\rho}(s)\\
%    &=\frac{\beta R}{n+\beta R}\rho_{0,n,m}-\varepsilon\frac{n}{n+\beta R}[\widehat{b}^1_{n,m}+\frac{\mu}{s+1}\widehat{b}^2_{n,m}+\beta\widehat{b}^3_{n,m}]\\
%    &+\varepsilon\frac{\mu}{s+1}\widehat{b}^2_{n,m}\frac{\sqrt{s+1}I_{n+3/2}(\sqrt{s+1}R)+\frac{n}{R} I_{n+1/2}(\sqrt{s+1}R)}{\sqrt{s+1}I_{n+3/2}(\sqrt{s+1}R)+(\frac{n}{R}+\beta) I_{n+1/2}(\sqrt{s+1}R)}\\
%    &-\frac{\partial}{\partial r}(\phi_{1,n,m}+\varepsilon\phi_{2,n,m})|_{r=R}+\frac{\mu}{s+1}(\frac{\partial\xi_{1,n,m}}{\partial r}+\varepsilon\frac{\partial\xi_{2,n,m}}{\partial r})\Big|_{r=R}+\varepsilon\widehat{b}^1_{n,m}\\
    &\frac{\beta R}{n+\beta R}\rho_{0,n,m}+\frac{\mu}{s+1}\frac{\partial\xi_{1,n,m}}{\partial r}\Big|_{r=R}-\frac{\partial\phi_{1,n,m}}{\partial r}\Big|_{r=R}+\varepsilon\Big\{\frac{\mu}{s+1}\frac{\partial\xi_{2,n,m}}{\partial r}\Big|_{r=R}-\frac{\partial\phi_{2,n,m}}{\partial r}\Big|_{r=R}\\
    &+\frac{\beta R}{n+\beta R}\widehat{b}^1_{n,m}+\frac{\beta R^2}{n+\beta R}\mu\widehat{b}^2_{n,m}\frac{1}{(s+1)R+(\frac{n}{R}+\beta)/P_n(\sqrt{s+1}R)}-\frac{n\beta }{n+\beta R}\widehat{b}^3_{n,m}\Big\},
\end{split}
\end{equation*}
where we have used the fact
\begin{equation*}
\begin{split}
  &\frac{\mu}{s+1}\widehat{b}^2_{n,m}\frac{\sqrt{s+1}I_{n+3/2}(\sqrt{s+1}R)+\frac{n}{R} I_{n+1/2}(\sqrt{s+1}R)}{\sqrt{s+1}I_{n+3/2}(\sqrt{s+1}R)+(\frac{n}{R}+\beta) I_{n+1/2}(\sqrt{s+1}R)}-\frac{\mu}{s+1}\widehat{b}^2_{n,m}\frac{n}{n+\beta R}\\
  %&=\frac{\mu}{s+1}\widehat{b}^2_{n,m}\Big(1-\frac{\beta I_{n+1/2}(\sqrt{s+1}R)}{\sqrt{s+1}I_{n+3/2}(\sqrt{s+1}R)+(\frac{n}{R}+\beta) I_{n+1/2}(\sqrt{s+1}R)}\Big)\\
%  &\quad-\frac{\mu}{s+1}\widehat{b}^2_{n,m}\Big(1-\frac{\beta R}{n+\beta R}\Big)\\
%  &=\frac{\mu}{s+1}\widehat{b}^2_{n,m}\frac{\beta R}{n+\beta R}\Big\{\frac{-\frac{n+\beta R}{R} I_{n+1/2}(\sqrt{s+1}R)}{\sqrt{s+1}I_{n+3/2}(\sqrt{s+1}R)+(\frac{n}{R}+\beta) I_{n+1/2}(\sqrt{s+1}R)}+1\Big\}\\
  &\quad=\frac{\mu}{s+1}\widehat{b}^2_{n,m}\frac{\beta R}{n+\beta R}\frac{\sqrt{s+1} I_{n+3/2}(\sqrt{s+1}R)}{\sqrt{s+1}I_{n+3/2}(\sqrt{s+1}R)+(\frac{n}{R}+\beta) I_{n+1/2}(\sqrt{s+1}R)}\\
  &\quad=\frac{\beta R^2}{n+\beta R}\mu\widehat{b}^2_{n,m}\frac{1}{(s+1)R+(\frac{n}{R}+\beta)/P_n(\sqrt{s+1}R)},\qquad (\mbox{by}\ (\ref{0.8})).
\end{split}
\end{equation*}
%As in \cite{HZH2}, the brace in left hand side of (\ref{0.53}) is equal to
%\begin{equation*}
%\begin{split}
%    &\frac{\beta R}{n+\beta R}\mu\frac{\beta RP_0(R)}{\beta+RP_0(R)}\left\{\frac{1}{\mu}\frac{\beta+RP_0(R)}{\beta RP_0(R)}\left[s+\frac{n}{R^3}\left(\frac{n(n+1)}{2}-1\right)\right]-RP_1(R)\right.\\
% &\qquad\qquad\qquad\left.+\frac{R^2P_1(R)+1+\beta R}{(s+1)R+(\frac{n}{R}+\beta)/P_n(\sqrt{s+1}R)}\right\}\\
%    &\equiv\frac{\beta R}{n+\beta R}\mu\frac{\beta RP_0(R)}{\beta+RP_0(R)}h_n(s,\mu,R),
%\end{split}
%\end{equation*}
It follows that
\begin{equation}\label{4.38}
\begin{split}
    \widehat{\rho}_{n,m}(s)
    %&\frac{n+\beta R}{\beta R}\frac{1}{\mu}\frac{\beta+RP_0(R)}{\beta RP_0(R)}\frac{1}{h_n(s,\mu,R)}\\
%    &\cdot\left\{\frac{\beta R}{n+\beta R}\rho_{0,n,m}-\varepsilon\frac{n}{n+\beta R}[\widehat{b}^1_{n,m}+\frac{\mu}{s+1}\widehat{b}^2_{n,m}+\beta\widehat{b}^3_{n,m}]\right.\\
%    &+\varepsilon\frac{\mu}{s+1}\widehat{b}^2_{n,m}\frac{\sqrt{s+1}I_{n+3/2}(\sqrt{s+1}R)+\frac{n}{R} I_{n+1/2}(\sqrt{s+1}R)}{\sqrt{s+1}I_{n+3/2}(\sqrt{s+1}R)+(\frac{n}{R}+\beta) I_{n+1/2}(\sqrt{s+1}R)}\\
%    &\left.-\frac{\partial}{\partial r}(\phi_{1,n,m}+\varepsilon\phi_{2,n,m})|_{r=R}+\frac{\mu}{s+1}(\frac{\partial\xi_{1,n,m}}{\partial r}+\varepsilon\frac{\partial\xi_{2,n,m}}{\partial r})\Big|_{r=R}+\varepsilon\widehat{b}^1_{n,m}\right\}\\
    =&\frac{n+\beta R}{\beta R}\frac{\beta+RP_0(R)}{\beta RP_0(R)}\frac{1}{\mu h_n(s,\mu,R)}\\
    &\cdot\Big\{\frac{\beta R}{n+\beta R}\rho_{0,n,m}+\frac{\mu}{s+1}\frac{\partial\xi_{1,n,m}}{\partial r}\Big|_{r=R}-\frac{\partial\phi_{1,n,m}}{\partial r}\Big|_{r=R}\\
    &\quad+\varepsilon\Big[\frac{\mu}{s+1}\frac{\partial\xi_{2,n,m}}{\partial r}\Big|_{r=R}-\frac{\partial\phi_{2,n,m}}{\partial r}\Big|_{r=R}+\frac{\beta R}{n+\beta R}\widehat{b}^1_{n,m}\\
  &\quad+\frac{\beta R^2}{n+\beta R}\mu\widehat{b}^2_{n,m}\frac{1}{(s+1)R+(\frac{n}{R}+\beta)/P_n(\sqrt{s+1}R)}-\frac{n\beta }{n+\beta R}\widehat{b}^3_{n,m}\Big]\Big\},
\end{split}
\end{equation}
which is equivalent to
\begin{equation}\label{4.71}
\begin{split}
  \widehat{\rho}_{n,m}(s)=&M_{n,m}(s)+\frac{n+\beta R}{\beta R}\frac{\beta+RP_0(R)}{\beta RP_0(R)}\frac{\varepsilon}{\mu h_n(s,\mu,R)}\\
  &\cdot\Big\{\frac{\mu}{s+1}\frac{\partial\xi_{2,n,m}}{\partial r}\Big|_{r=R}-\frac{\partial\phi_{2,n,m}}{\partial r}\Big|_{r=R}+\frac{\beta R}{n+\beta R}\widehat{b}^1_{n,m}\\
  &\quad+\frac{\beta R^2}{n+\beta R}\mu\widehat{b}^2_{n,m}\frac{1}{(s+1)R+(\frac{n}{R}+\beta)/P_n(\sqrt{s+1}R)}-\frac{n\beta }{n+\beta R}\widehat{b}^3_{n,m}\Big\},
\end{split}
\end{equation}
where
\begin{equation}\label{0.60}
\begin{split}
  M_{n,m}(s)=&\frac{n+\beta R}{\beta R}\frac{\beta+RP_0(R)}{\beta RP_0(R)}\frac{1}{\mu h_n(s,\mu,R)}\\
  &\cdot\Big\{\frac{\beta R}{n+\beta R}\rho_{0,n,m}+\frac{\mu}{s+1}\frac{\partial\xi_{1,n,m}}{\partial r}\Big|_{r=R}-\frac{\partial\phi_{1,n,m}}{\partial r}\Big|_{r=R}\Big\}.
\end{split}
\end{equation}

\section{A choice of a new center}
After the initial values perturbed, say (\ref{1.13}), the domain of the system (\ref{3.2})-(\ref{3.5}) (or (\ref{1.2})-(\ref{1.6})) undergoes a translation owing to the perturbation of ``hidden" mode 1 contributions. In order to prove asymptotic stability we need to determine the center of the limiting sphere.
Unlike the linear stability, it is a grand challenge to find the new center. In this section, we shall accomplish this through the following contraction mapping type of theorem:
\begin{thm}\label{th2}
Let $(X,\|\cdot\|)$ be a Banach space and let $\overline{B}_K(a_0)$ denote the closed ball in $X$ with center $a_0$ and radius $K$. Let $F$ be a mapping from $\overline{B}_K(a_0)$ into $X$ and
\begin{equation*}
F(x)=F_1(x)+\varepsilon G(x),
\end{equation*}
such that

(i) $F_1'(x)$ and $G'(x)$ are both continuous for $x\in\overline{B}_K(a_0)$,

(ii) $F_1(a_0)=0$ and the operator $F'_1(a_0)$ is invertible.

\noindent Then for small $|\varepsilon|$, the equation $F(x)=0$ admits a unique solution $x$ in $\overline{B}_K(a_0)$.
\end{thm}

\noindent{\bf Proof.}\  $F(x)=0$ if and only if $F_1(x)=-\varepsilon G(x)$. By implicit function theorem, for $\varepsilon$ small, this is equivalent to
\begin{equation*}
  x=F_1^{-1}(-\varepsilon G(x)).
\end{equation*}
Let $g(x)=F_1^{-1}(-\varepsilon G(x))$, then
\begin{equation*}
  \|g'(x)\|\leq\|(F^{-1}_1)'(-\varepsilon G(x))\|\cdot|\varepsilon|\cdot\|G'(x)\|.
\end{equation*}
Under our assumptions, $\|G'(x)\|\leq C$ for $x\in\overline{B}_K(a_0)$ and $(F_1^{-1})'(y)$ is continuous for $y$ in a small neighborhood of 0, therefore
\begin{equation*}
  \|(F_1^{-1})'(-\varepsilon G(x))\|\leq C, \ \ \ \ \ x\in\overline{B}_K(a_0)
\end{equation*}
for $|\varepsilon|$ sufficiently small. It follows that $\|g'(x)\|\leq C|\varepsilon|$ for $x\in\overline{B}_K(a_0)$. Taking $|\varepsilon|$ to be small enough, we get a contraction, which provides a unique fixed point for $g(x)$ in $\overline{B}_K(a_0)$. \hfill $\Box$

For $n\neq1$, all zeros of $h_n(s)$ lie in $\mbox{Re}\ s<-\delta(n+1)$ with some $\delta>0$ (see \cite{HZH2}), which enables us to change the contour $\Gamma$ of integration to $-\delta(n+1)$ in the inverse Laplace transform of the various functions (\ref{4.71}). However, for $n=1$, since $h_1(s)$ has a simple zero at $s=0$ (see also \cite{HZH2}), we cannot move the contour $\Gamma$ to $\mbox{Re}\ s=-\delta$ with some $\delta>0$ in the inverse Laplace transform and obtain the decay rate in $t$. In the sequel, we will make use of Theorem \ref{th2} to show that there is a translation of coordinates
\begin{equation*}
  0\rightarrow\varepsilon a^\ast(\varepsilon),
\end{equation*}
where $a^\ast(\varepsilon)$ is uniformly bounded, such that for $n=1$, the expression in the brace in (\ref{4.38}) vanishes at $s=0$ and thus the singularity of $1/h_1(s)$ at $s=0$ will be cancelled. Hence, in the new coordinate system, we can take the inverse Laplace transform in (\ref{4.38}) and move the contour $\Gamma$ to $\mbox{Re}\ s=-\delta$.

In the process of using the contraction mapping (Theorem \ref{th2}), we adjust the center of the domain to an optimal location $a$ at each iteration. For notational convenience, we use the same letter $r$, $\theta$, $\varphi$ to denote the new variables after the translation of coordinates $0\rightarrow\varepsilon a$ with $a=(a_1,a_2,a_3)$. As in the proof of \cite[(6.4) and (6.5) in Lemma 6.1]{FH2}, the initial data $w_0$, $\rho_0$ in the new coordinates are transformed into, respectively,
\begin{equation}\label{0.62}
  w_0(r,\theta,\varphi)+\frac{\partial\sigma_S(r)}{\partial r}(a_1\cos\varphi\sin\theta+a_2\sin\varphi\sin\theta+a_3\cos\theta)+\varepsilon A(r,\theta,\varphi,\varepsilon;a),
\end{equation}
\begin{equation}\label{0.63}
  \rho_0(\theta,\varphi)-(a_1\cos\varphi\sin\theta+a_2\sin\varphi\sin\theta+a_3\cos\theta)+\varepsilon B(r,\theta,\varphi,\varepsilon;a),
\end{equation}
where $A$ and $B$ are bounded functions. Recall that
\begin{equation*}
  Y_{1,-1}=\sqrt{\frac{3}{8\pi}}\sin\theta e^{-i\varphi},\ \ \ \ Y_{1,0}=\sqrt{\frac{3}{4\pi}}\cos\theta,\ \ \ \ Y_{1,1}=-\sqrt{\frac{3}{8\pi}}\sin\theta e^{i\varphi},
\end{equation*}
 and $w_{0,1,m}$, $\rho_{0,1,m}$ $(m=-1,0,1)$ are changed into (see \cite[(6.6) and (6.7)]{FH2})
\begin{equation}\label{6.2}
    \widetilde{w}_{0,1,m}\triangleq w_{0,1,m}+b_{m+2}\frac{\partial\sigma_S(r)}{\partial r}+\varepsilon A_m,
      %\ \ \ \ \mbox{and}\ \ \ \ \ \widetilde{\rho}_{0,1,m}\equiv\rho_{0,1,m}-b_{m+2}+\varepsilon B_m,
\end{equation}
\begin{equation}\label{0.66}
  \widetilde{\rho}_{0,1,m}\triangleq\rho_{0,1,m}-b_{m+2}+\varepsilon B_m,
\end{equation}
where $A_m$ and $B_m$ are bounded functions of $(a,\varepsilon)$, and $b_{m+2}$ satisfies (see \cite[(6.8)]{FH2})
\begin{equation*}
 b_1-b_3=a_1\sqrt{\frac{8\pi}{3}}, \ \ \ \ i(b_1+b_3)=-a_2\sqrt{\frac{8\pi}{3}}, \ \ \ \ b_2=a_3\sqrt{\frac{4\pi}{3}}.
\end{equation*}
%Note that after the initial values perturbed, say (\ref{1.13}), the system (\ref{3.2})-(\ref{3.6}) (or (\ref{1.2})-(\ref{1.6})) undergoes a translation of the origin resulting from the perturbation of mode 1. In order to prove asymptotic stability, we need to find the new center, which will be discussed in Section 4. We also note that the translation of the origin does not change the equation (\ref{3.2})-(\ref{3.6}).
Note that a translation of the origin does not change the equations (\ref{3.2})-(\ref{3.5}), but change the initial data. Namely, moving the origin of the system (\ref{3.2})-(\ref{3.5}) to $\varepsilon a$ is equivalent to keep the origin fixed at $0$ but replace the initial values $w_0$ and $\rho_0$ by (\ref{0.62}) and (\ref{0.63}), respectively. %Furthermore, as we go to (\ref{3.15})-(\ref{3.18}), we shall always use the translated initial data.

As stated in \cite[page 632]{FH2}, we note that the functions $(f^i,b^j)$ satisfying (\ref{4.1})-(\ref{4.10}) depend on the new center $a$, i.e., $(f^i(r,\theta,\varphi,t;a),b^j(\theta,\varphi,t;a))$, to ensure the consistency condition of order 2.

In the new coordinate system, for $n=1$, $s=0$, the expression in the braces in (\ref{4.38}) is given by
\begin{equation}\label{3.30}
\begin{split}
  F_m(a)\equiv\Big\{&\frac{\beta R}{1+\beta R}\widetilde{\rho}_{0,1,m}+\frac{\mu}{s+1}\frac{\partial\xi_{1,1,m}}{\partial r}\Big|_{r=R}-\frac{\partial\phi_{1,1,m}}{\partial r}\Big|_{r=R}\\
  &+\varepsilon\Big[\frac{\mu}{s+1}\frac{\partial\xi_{2,1,m}}{\partial r}\Big|_{r=R}-\frac{\partial\phi_{2,1,m}}{\partial r}\Big|_{r=R}+\frac{\beta R}{1+\beta R}\widehat{b}^1_{1,m}\\
  &+\frac{\beta R^2}{1+\beta R}\mu\widehat{b}^2_{1,m}\frac{1}{(s+1)R+(\frac{1}{R}+\beta)/P_1(\sqrt{s+1}R)}-\frac{\beta }{1+\beta R}\widehat{b}^3_{1,m}\Big]\Big\}\Big|_{s=0},
\end{split}
\end{equation}
where $m=-1,0,1$. In fact, $F_m(a)$ depends implicitly on $a$ through the dependence of each term of the right-hand side on $a$.
%Let $F(a)=(F_{-1}(a),F_0(a),F_1(a))$.

We first establish the following lemma.
\begin{lem}\label{l6.1}
%Let
%\begin{equation*}
%\begin{split}
%  (f^1,f^2)=(f^1(r,\theta,\varphi,t;a),f^2(r,\theta,\varphi,t;a)),
%\end{split}
%\end{equation*}
%\begin{equation*}
%  (b^1,b^2,b^3)=(b^1(\theta,\varphi,t;a),b^2(\theta,\varphi,t;a),b^3(\theta,\varphi,t;a))
%\end{equation*}
%be any functions satisfying (\ref{4.1})-(\ref{4.10}).
Let $\mu<\mu^\ast$, then there exists a translation of the origin $0\rightarrow\varepsilon a$ with $a=(a_1,a_2,a_3)$ uniformly bounded, i.e. $|a|\leq Q_0$, such that in the new coordinate system,
\begin{equation}\label{6.1}
  \left\{\frac{\beta R}{1+\beta R}\widetilde{\rho}_{0,1,m}+\frac{\mu}{s+1}\frac{\partial\xi_{1,1,m}}{\partial r}\Big|_{r=R}-\frac{\partial\phi_{1,1,m}}{\partial r}\Big|_{r=R}\right\}\Big|_{s=0}=\varepsilon C_m, \ \ \ \ m=-1,0,1,
\end{equation}
where $C_m$ is a bounded function of $(a,\varepsilon)$.
\end{lem}
\noindent{\bf Proof.}\
%As in the proof of \cite[Lemma 6.1]{FH2}, in the new coordinate system the initial data $w_0$ and $\rho_0$ are changed into
%\begin{equation}\label{6.2}
%  \widetilde{\rho}_{0,1,m}\equiv\rho_{0,1,m}-b_{m+2}+\varepsilon A_m, \ \ \ \ \mbox{and}\ \ \ \ \ \widetilde{w}_{0,1,m}\equiv w_{0,1,m}+b_{m+2}\frac{\partial\sigma_S(r)}{\partial r}+\varepsilon B_m,
%\end{equation}
%where $A_m$ and $B_m$ are bounded functions of $(a,\varepsilon)$, and $b_{m+2}$ depends on $a$.
By (\ref{4.30}), (\ref{4.35}) and (\ref{6.2}), in the new coordinate system, the functions $\xi_{1,1,m}$ and $\phi_{1,1,m}$ at $s=0$ satisfy
\begin{equation}\label{6.3}
\begin{split}
  -\Delta\xi_{1,1,m}+\left(\frac{2}{r^2}+1\right)\xi_{1,1,m} &= w_{0,1,m}(r)+b_{m+2}\frac{\partial\sigma_S(r)}{\partial r}+\varepsilon A_m \ \ \ \ \mbox{in}\ B_R, \\
  \left(\frac{\partial \xi_{1,1,m}}{\partial r}+\beta \xi_{1,1,m}\right)\Big|_{r=R} &= 0,
\end{split}
\end{equation}
\begin{equation}\label{6.4}
\begin{split}
  -\Delta\phi_{1,1,m}+\frac{2}{r^2}\phi_{1,1,m} &= \mu\Big(w_{0,1,m}(r)+b_{m+2}\frac{\partial\sigma_S(r)}{\partial r}+\varepsilon A_m\Big)\ \ \ \ \mbox{in}\ B_R, \\
  \left(\frac{\partial\phi_{1,1,m}}{\partial r}+\beta\phi_{1,1,m}\right)\Big|_{r=R} &= 0.
\end{split}
\end{equation}
The functions $\xi_{1,1,m}$ and $\phi_{1,1,m}$ can be rewritten as
\begin{equation*}
  \xi_{1,1,m}=b_{m+2}\xi_{1,1,m}^{(b)}+\widetilde{\xi}_{1,1,m}, \ \ \ \ \phi_{1,1,m}=b_{m+2}\phi_{1,1,m}^{(b)}+\widetilde{\phi}_{1,1,m},
\end{equation*}
where $\widetilde{\xi}_{1,1,m}$, $\widetilde{\phi}_{1,1,m}$ are the solutions of (\ref{6.3}) and (\ref{6.4}) corresponding to the right-hand side terms $w_{0,1,m}+\varepsilon A_m$ and $\mu(w_{0,1,m}+\varepsilon A_m)$, respectively.

Substituting these results and (\ref{0.66}) into (\ref{6.1}), we obtain
\begin{equation}\label{3.31}
\begin{split}
  &\Big\{\frac{\beta R}{1+\beta R}\widetilde{\rho}_{0,1,m}+\frac{\mu}{s+1}\frac{\partial\xi_{1,1,m}}{\partial r}\Big|_{r=R}-\frac{\partial\phi_{1,1,m}}{\partial r}\Big|_{r=R}\Big\}\Big|_{s=0}\\
  &=\Big\{-\frac{\beta R}{1+\beta R}+\mu\frac{\partial\xi^{(b)}_{1,1,m}}{\partial r}\Big|_{r=R}-\frac{\partial\phi^{(b)}_{1,1,m}}{\partial r}\Big|_{r=R}\Big\}b_{m+2}+G_m^1+\varepsilon G_m^2\\
  &\equiv Q_mb_{m+2}+G_m^1+\varepsilon G_m^2,
\end{split}
\end{equation}
where $G_m^1$ is independent of $b_{m+2}$ and $G_m^2$ is a bounded function of $a$ and $\varepsilon$. If we can show $Q_m\neq0$, then we can choose $b_{m+2}$ to cancel $G_m^1$, leaving out the expression $\varepsilon G_m^2$, which is the right-hand side of (\ref{6.1}). In fact, we have the following lemma:
\begin{lem}\label{l0.1}
$Q_m=0$ if and only if $\mu=\mu_1(R)$, where $\mu_1(R)$ is defined by (\ref{0.035}).
\end{lem}
\noindent{\bf Proof.}\ In order to compute $Q_m$, we need to compute the functions $\xi_{1,1,m}^{(b)}$ and $\phi_{1,1,m}^{(b)}$ which satisfy
\begin{equation}\label{6.5}
\begin{split}
  -\Delta\xi_{1,1,m}^{(b)}+\Big(\frac{2}{r^2}+1\Big)\xi_{1,1,m}^{(b)} = \frac{\partial\sigma_S(r)}{\partial r}\ \ \ \ \mbox{in}\ B_R,\ \ \ \Big(\frac{\partial \xi_{1,1,m}^{(b)}}{\partial r}+\beta \xi_{1,1,m}^{(b)}\Big)\Big|_{r=R} = 0,
\end{split}
\end{equation}
%\begin{equation}\label{6.6}
%  \Big(\frac{\partial \xi_{1,1,m}^{(b)}}{\partial r}+\beta \xi_{1,1,m}^{(b)}\Big)\Big|_{r=R} = 0,
%\end{equation}
and
\begin{equation}\label{6.7}
\begin{split}
  -\Delta\phi_{1,1,m}^{(b)}+\frac{2}{r^2}\phi_{1,1,m}^{(b)} = \mu\frac{\partial\sigma_S(r)}{\partial r} \ \ \mbox{in}\ B_R,
  \ \ \ \ \Big(\frac{\partial \phi_{1,1,m}^{(b)}}{\partial r}+\beta \phi_{1,1,m}^{(b)}\Big)\Big|_{r=R} = 0.
\end{split}
\end{equation}
%\begin{equation}\label{6.8}
%  \Big(\frac{\partial \phi_{1,1,m}^{(b)}}{\partial r}+\beta \phi_{1,1,m}^{(b)}\Big)\Big|_{r=R} = 0.
%\end{equation}

Using the equation for $\sigma'_S(r)$, i.e. $\sigma'''_S+\frac{2}{r}\sigma''_S-\frac{2}{r^2}\sigma'_S=\sigma'_S$, one can easily verify that
\begin{equation*}
  \phi_{1,1,m}^{(b)}=-\mu\sigma'_S(r)+\frac{\mu\lambda}{1+\beta R}r,\ \ \ \ \ (\lambda=\sigma_S''(R)+\beta\sigma_S'(R)\ \mbox{by}\ (\ref{0.22})),
\end{equation*}
then by (\ref{1.7}) and (\ref{0.22}), we have
\begin{equation}\label{6.9}
\begin{split}
  %\frac{\partial\phi_{1,1,m}^{(b)}(R)}{\partial r}=-\beta\phi_{1,1,m}^{(b)}(R)&=\mu\beta\sigma'_S(R)-\frac{\mu\lambda}{1+\beta R}\beta R\\
  %&=\mu\beta\frac{\beta RP_0(R)}{\beta+RP_0(R)}-\frac{\mu\lambda}{1+\beta R}\beta R.
  \frac{\partial\phi_{1,1,m}^{(b)}(R)}{\partial r}&=-\mu\sigma''_S(R)+\frac{\mu\lambda}{1+\beta R}\\
  &=-\mu\frac{\beta[1-2P_0(R)]}{\beta+RP_0(R)}+\frac{\mu}{1+\beta R}\frac{\beta P_0(R)}{\beta+RP_0(R)}[R^2P_1(R)+1+\beta R].
\end{split}
\end{equation}

We next proceed to compute $\xi_{1,1,m}^{(b)}$. As in \cite[page 630]{FH2}, the equation for $\xi_{1,1,m}^{(b)}$ can be rewritten in the following form:
\begin{equation*}
  -r^2[\sigma_S'(r)]^2\frac{\partial}{\partial r}\Big(\frac{\xi_{1,1,m}^{(b)}}{\sigma_S'(r)}\Big)=r^2\sigma_S(r)\sigma_S'(r)-\int_0^r\tau^2\sigma_S^2(\tau)d\tau,
\end{equation*}
then it follows that
\begin{equation}\label{6.10}
  \frac{\partial_r\xi_{1,1,m}^{(b)}\sigma_S'(r)-\xi_{1,1,m}^{(b)}\sigma_S''(r)}{[\sigma_S'(r)]^2}=-\frac{\sigma_S(r)}{\sigma_S'(r)}
  +\frac{1}{r^2[\sigma_S'(r)]^2}\int_0^r\tau^2\sigma_S^2(\tau)d\tau.
\end{equation}
By (\ref{6.10}), the boundary condition of $\xi^{(b)}_{1,1,m}$ in (\ref{6.5}) and (\ref{0.22}), we obtain
\begin{equation}\label{6.11}
  \frac{\partial\xi_{1,1,m}^{(b)}(R)}{\partial r}=\frac{\beta+RP_0(R)}{P_0(R)[R^2P_1(R)+1+\beta R]}\Big[-\sigma_S(R)\sigma'_S(R)+\frac{1}{R^2}\int_0^R\tau^2\sigma_S^2(\tau)d\tau\Big].
\end{equation}
%Using (\ref{6.9}) and (\ref{6.11}), we compute
%\begin{equation}\label{6.12}
%\begin{split}
%  Q_m=&-\frac{\beta R}{1+\beta R}+\frac{\partial}{\partial r}\Big[\frac{\mu}{s+1}\xi^{(b)}_{1,1,m}-\phi^{(b)}_{1,1,m}\Big]\Big|_{r=R}\\
%  =&-\frac{\beta R}{1+\beta R}+\frac{\mu\beta}{\lambda}[-\sigma_S(R)\sigma'_S(R)+\frac{1}{R^2}\int_0^R\tau^2(\sigma_S)^2(\tau)d\tau]\\
%  &+\mu[\sigma_S(R)-\frac{2}{R}\sigma'_S(R)]-\frac{\mu\lambda}{1+\beta R}.
%\end{split}
%\end{equation}
It follows from (\ref{1.7}) and (\ref{0.67}) that the second term of the right-hand side of (\ref{6.11}) is equal to
\begin{equation*}
\begin{split}
  \int_0^R\tau^2\sigma_S^2(\tau)d\tau=&\Big(\frac{\beta}{\beta+RP_0(R)}\frac{R^{1/2}}{I_{1/2}(R)}\Big)^2\int_0^R\tau I_{1/2}^2(\tau)d\tau\\
  =&\Big(\frac{\beta}{\beta+RP_0(R)}\frac{R^{1/2}}{I_{1/2}(R)}\Big)^2\int_0^R\frac{2}{\pi}\sinh^2\tau d\tau\\
  %=&\Big(\frac{\beta}{\beta+RP_0(R)}\frac{R^{1/2}}{I_{1/2}(R)}\Big)^2\frac{2}{\pi}\frac{\sinh^2R}{2}\Big(\coth R-\frac{R}{\sinh^2R}\Big)   \\
  %=&\Big(\frac{\beta}{\beta+RP_0(R)}\frac{R^{1/2}}{I_{1/2}(R)}\Big)^2\frac{2}{\pi}\frac{\sinh^2R}{2}(\coth R-R\coth^2R+R)\\
  =&\Big(\frac{\beta}{\beta+RP_0(R)}\frac{R^{1/2}}{I_{1/2}(R)}\Big)^2\frac{2}{\pi}\frac{\sinh^2R}{2}R\big[-P_0(R)-R^2P_0^2(R)+1\big]\\
  =&\Big(\frac{\beta}{\beta+RP_0(R)}\Big)^2\frac{R^3}{2}\big[-P_0(R)-R^2P_0^2(R)+1\big],
\end{split}
\end{equation*}
where $I_{1/2}^2(\tau)=\frac{2}{\pi \tau}\sinh^2\tau$ is used, so that, by (\ref{1.7}) and (\ref{6.14}),
\begin{equation}\label{6.13}
\begin{split}
  \frac{\partial\xi_{1,1,m}^{(b)}(R)}{\partial r}=&\frac{\beta+RP_0(R)}{P_0(R)[R^2P_1(R)+1+\beta R]}\Big(\frac{\beta}{\beta+RP_0(R)}\Big)^2\\
  &\cdot\Big[-\frac{3}{2}RP_0(R)-\frac{1}{2}R^3P_0^2(R)+\frac{R}{2}\Big].
\end{split}
\end{equation}

Substituting the above results (\ref{6.9}) and (\ref{6.13}) into the expression $Q_m$, we derive
%\begin{equation*}
%\begin{split}
%  Q_m&=-\frac{\beta R}{1+\beta R}+\frac{\partial}{\partial r}\Big[\mu\xi^{(b)}_{1,1,m}-\phi^{(b)}_{1,1,m}\Big]\Big|_{r=R}\\
%  &=-\frac{\beta R}{1+\beta R}+\mu \frac{\beta}{\lambda}\Big(\frac{\beta}{\beta+RP_0(R)}\Big)^2\Big[-\frac{3}{2}RP_0(R)-\frac{1}{2}R^3P_0^2(R)+\frac{R}{2}\Big]\\
%  &\quad-\mu\Big[\frac{\beta^2RP_0(R)}{\beta+RP_0(R)}-\frac{\lambda}{1+\beta R}\beta R\Big]\\
%  &=-\frac{\beta R}{1+\beta R}+\mu\Big\{\frac{\beta}{P_0(R)(R^2P_1(R)+1+\beta R)}\frac{\beta}{\beta+RP_0(R)}\\
%  &\quad\cdot\Big[-\frac{3}{2}RP_0(R)-\frac{1}{2}R^3P_0^2(R)+\frac{R}{2}\Big]-\beta\frac{\beta RP_0(R)}{\beta+RP_0(R)}\\
%  &\quad+\frac{\beta RP_0(R)}{\beta+RP_0(R)}(R^2P_1(R)+1+\beta R)\frac{\beta}{1+\beta R}\Big\}\ \ \ \ \ \ \ \  (\mbox{by}\ (\ref{0.22}))\\
%  %&=-\frac{\beta R}{1+\beta R}+\mu\frac{\beta R}{\beta+RP_0(R)}\{\frac{\beta}{R^2P_1(R)+1+\beta R}[-\frac{3}{2}-\frac{1}{2}R^2P_0(R)+\frac{1}{2P_0(R)}]\\
%%  &-\beta P_0(R)+(R^2P_1(R)+1+\beta R)\frac{\beta P_0(R)}{1+\beta R}\}\\
%  &=-\frac{\beta R}{1+\beta R}+\mu\frac{\beta R}{\beta+RP_0(R)}\Big\{\frac{\beta}{R^2P_1(R)+1+\beta R}\Big[-\frac{3}{2}-\frac{1}{2}R^2P_0(R)+\frac{1}{2P_0(R)}\Big]\\
%  &\quad+\frac{\beta R^2P_0(R)P_1(R)}{1+\beta R}\Big\}\\
%  &=\frac{\beta R}{1+\beta R}\Big\{-1+\frac{\mu\beta}{\beta+RP_0(R)}\\
%  &\quad\cdot\frac{(1+\beta R)(-\frac{3}{2}-\frac{1}{2}R^2P_0(R)+\frac{1}{2P_0(R)})+ R^2P_0(R)P_1(R)(R^2P_1(R)+1+\beta R)}{R^2P_1(R)+1+\beta R}\Big\}.
%\end{split}
%\end{equation*}
\begin{equation*}
\begin{split}
  Q_m%&=-\frac{\beta R}{1+\beta R}+\mu\frac{\partial\xi^{(b)}_{1,1,m}}{\partial r}\Big|_{r=R}-\frac{\partial\phi^{(b)}_{1,1,m}}{\partial r}\Big|_{r=R}\\
  &=-\frac{\beta R}{1+\beta R}\\
  &\quad+\mu\frac{\beta+RP_0(R)}{P_0(R)[R^2P_1(R)+1+\beta R]}\Big(\frac{\beta}{\beta+RP_0(R)}\Big)^2\Big[-\frac{3}{2}RP_0(R)-\frac{1}{2}R^3P_0^2(R)+\frac{R}{2}\Big]\\
  &\quad+\mu\frac{\beta}{\beta+RP_0(R)}[1-2P_0(R)]
  -\frac{\mu}{1+\beta R}\frac{\beta P_0(R)}{\beta+RP_0(R)}[R^2P_1(R)+1+\beta R]\\
  %&=-\frac{\beta R}{1+\beta R}+\mu\frac{\beta}{\beta+RP_0(R)}\\
%  &\quad\cdot\Big\{\frac{\beta}{P_0(R)[R^2P_1(R)+1+\beta R]}\Big[-\frac{3}{2}RP_0(R)-\frac{1}{2}R^3P_0^2(R)+\frac{R}{2}\Big]\\
%  &\qquad+1-2P_0(R)-\frac{P_0(R)}{1+\beta R}[R^2P_1(R)+1+\beta R]\Big\}\\
  &=-\frac{\beta R}{1+\beta R}+\mu\frac{\beta R}{1+\beta R}\frac{\beta}{\beta+RP_0(R)}\\
  &\quad\cdot\Big\{\frac{1+\beta R}{R^2P_1(R)+1+\beta R}\Big[-\frac{3}{2}-\frac{1}{2}R^2P_0(R)+\frac{1}{2P_0(R)}\Big]\\
  &\quad+\frac{1+\beta R}{\beta R}[1-2P_0(R)]-\frac{R^2P_0(R)P_1(R)}{\beta R}-\frac{1+\beta R}{\beta R}P_0(R)\Big\}\\
  %&=-\frac{\beta R}{1+\beta R}+\mu\frac{\beta R}{1+\beta R}\frac{\beta}{\beta+RP_0(R)}\\
%  &\quad\cdot\Big\{\frac{1+\beta R}{R^2P_1(R)+1+\beta R}\Big[-\frac{3}{2}-\frac{1}{2}R^2P_0(R)+\frac{1}{2P_0(R)}\Big]\\
%  &\qquad+\frac{1+\beta R}{\beta R}[1-3P_0(R)]-\frac{R^2P_0(R)P_1(R)}{\beta R}\Big\}\\
  &=-\frac{\beta R}{1+\beta R}+\mu\frac{\beta R}{1+\beta R}\frac{\beta}{\beta+RP_0(R)}\\
  &\quad\cdot\Big\{\frac{1+\beta R}{R^2P_1(R)+1+\beta R}\Big[-\frac{3}{2}-\frac{1}{2}R^2P_0(R)+\frac{1}{2P_0(R)}\Big]+R^2P_0(R)P_1(R)\Big\}\ \ (\mbox{by}\ (\ref{0.10}))\\
  &=\frac{\beta R}{1+\beta R}\Big\{-1+\mu\frac{\beta}{\beta+RP_0(R)}\frac{1}{R^2P_1(R)+1+\beta R}\\
  &\quad\cdot\Big[(1+\beta R)\Big(-\frac{3}{2}-\frac{1}{2}R^2P_0(R)+\frac{1}{2P_0(R)}\Big)+ R^2P_0(R)P_1(R)(R^2P_1(R)+1+\beta R)\Big]\Big\}
\end{split}
\end{equation*}

Furthermore, the bracket in the right-hand side can simplify to
\begin{equation*}
\begin{split}
  (1+\beta R)&\Big(-\frac{3}{2}-\frac{1}{2}R^2P_0(R)+\frac{1}{2P_0(R)}\Big)+ R^2P_0(R)P_1(R)(R^2P_1(R)+1+\beta R)\\
  =&(1+\beta R)\Big(-\frac{3}{2}-\frac{1}{2}R^2P_0(R)+\frac{1}{2P_0(R)}\Big)+R^4P_0(R)P_1^2(R)\\
  &+(1+\beta R)R^2P_0(R)P_1(R)\\
  =&(1+\beta R)\Big(-\frac{3}{2}-\frac{1}{2}R^2P_0(R)+\frac{1}{2P_0(R)}\Big)+R^4P_0(R)P_1^2(R)\\
  &+(1+\beta R)(1-3P_0(R)) \qquad \qquad \qquad \ \ \ (\mbox{by}\ (\ref{0.10}))\\
  =&(1+\beta R)\Big(-\frac{R^2P_0(R)}{2}-\frac{1}{2}-3P_0(R)+\frac{1}{2P_0(R)}\Big)+R^4P_0(R)P_1^2(R)\\
  =&-\frac{1+\beta R}{2}\Big[R^2P_0(R)+1+6P_0(R)-\frac{1}{P_0(R)}\Big]+R^4P_0(R)P_1^2(R)\\
  =&-\frac{1+\beta R}{2}R^3P_0(R)P_1'(R)+R^4P_0(R)P_1^2(R)\\
  =&R^3P_0(R)\Big[RP_1^2(R)-\frac{1+\beta R}{2}P_1'(R)\Big],
\end{split}
\end{equation*}
where we have used the fact (see \cite[page 631]{FH2})
\begin{equation*}
  R^2P_0(R)+1+6P_0(R)-\frac{1}{P_0(R)}=R^3P_0(R)P_1'(R).
\end{equation*}
Hence, by (\ref{0.035}), we obtain
\begin{equation}\label{0.3}
\begin{split}
  Q_m&=\frac{\beta R}{1+\beta R}\Big\{-1+\mu\frac{\beta R^3P_0(R)}{\beta+RP_0(R)}\cdot\frac{RP_1^2(R)-\frac{1+\beta R}{2}P_1'(R)}{R^2P_1(R)+1+\beta R}\Big\}\\
  &=\frac{\beta R}{1+\beta R}\Big\{-1+\frac{\mu}{\mu_1(R)}\Big\},
\end{split}
\end{equation}
so that $Q_m=0$ if and only if $\mu=\mu_1(R)$,
%\begin{equation*}
%  Q_m=0 \ \ \ \mbox{if\ and\ only\ if} \ \ \ \ \ \mu=\mu_1(R),
%\end{equation*}
which completes the proof of this lemma.
\hfill $\Box$

Hence, recalling (\ref{0.036}), i.e. $\mu<\mu^\ast<\mu_1(R)$, we complete the proof of Lemma \ref{l6.1}.   \hfill $\Box$

%Note that a translation of coordinates does not change the equation (\ref{3.2})-(\ref{3.6}), but change the initial values.

%Let $F(a)=(F_{-1}(a),F_0(a),F_1(a))$, where
%\begin{equation}\label{3.32}
%\begin{split}
%  F_m(a)=&\left\{\frac{\beta R}{1+\beta R}\widetilde{\rho}_{0,1,m}+\frac{\mu}{s+1}\frac{\partial\xi_{1,1,m}}{\partial r}|_{r=R}-\frac{\partial\phi_{1,1,m}}{\partial r}|_{r=R}\right.\\
%  &\quad+\varepsilon\Big[\frac{\mu}{s+1}\frac{\partial\xi_{2,1,m}}{\partial r}|_{r=R}-\frac{\partial\phi_{2,1,m}}{\partial r}|_{r=R}\Big]\\
%  &\quad+\varepsilon\Big[\frac{\beta R}{1+\beta R}\widehat{b}^1_{1,m}+\frac{\beta R^2}{1+\beta R}\mu\widehat{b}^2_{1,m}\frac{1}{(s+1)R+(\frac{1}{R}+\beta)/P_1(\sqrt{s+1}R)}\\
%  &\qquad\quad\left.-\frac{\beta }{1+\beta R}\widehat{b}^3_{1,m}\Big]\right\}\Big|_{s=0}, \ \ \ \ m=-1,0,1.
%\end{split}
%\end{equation}
%In fact, $F(a)$ depends implicitly on $a$ through the dependence of each term of the right-hand side on $a$.
%$\widetilde{\rho}_{0,1,m}$, $\xi_{1,1,m}$, $\phi_{1,1,m}$, $\xi_{2,1,m}$, $\phi_{2,1,m}$ and $\widehat{b}^j_{1,m}$ on $a$.

Let $F(a)=(F_{-1}(a),F_0(a),F_1(a))$. Now we proceed to give the main result of this section.
\begin{thm}\label{th1}
There exists a new center $\varepsilon a^\ast(\varepsilon)$ such that $F(a^\ast(\varepsilon))=0$.
\end{thm}
\noindent{\bf Proof.}\ By (\ref{3.30}) and (\ref{3.31}),
\begin{equation*}
F(a)=E(a)+\varepsilon G(a).
\end{equation*}
It follows from Lemma \ref{l6.1} that we can choose $a_0$ such that $E(a_0)=0$. Clearly, $E'(a)$ and $G'(a)$ are continuous. Since $Q_m\neq0$, $E'(a_0)$ is invertible. Therefore, by Theorem \ref{th2}, the proof is complete.
 \hfill $\Box$

\section{The Inhomogeneous Linear System}
In this section, we shall derive the estimates of $\rho(\theta,\varphi,t)$ for the inhomogeneous linear system (\ref{3.15})-(\ref{3.18}) which is in the new coordinate system.
%As in Section 3, by spherical extension %$\rho(\theta,\varphi,t)=\sum_{n=0}^\infty\sum_{m=-n}^n\rho_{n,m}(t)Y_{n,m}(\theta,\varphi)$
%and Laplace transforms,
%$\widehat{\rho}_{n,m}(s)=\int^{\infty}_0e^{-st}\rho_{n,m}(t)dt$,
%the system (\ref{3.15})-(\ref{3.18}) can be solved explicitly and $\widehat{\rho}_{n,m}(s)$ has the form (\ref{4.38}).
In order to do so, we need to take the inverse Laplace transform of $\widehat{\rho}_{n,m}(s)$ which is of the form (\ref{4.38}) (or (\ref{4.71})). Since for $n=1$, $h_1(s)$ has a simple root at $s=0$ that is different from the situation $n\neq1$, the computation is divided into two cases: $n\neq1$ and $n=1$.

\subsection{Case 1: $n\neq1$}

%As a matter of fact, the Laplace transform of $\rho_{n,m}(t)$ has the representation (\ref{4.71}).
In this case, the choice of new center does not change the proof.
Introduce the following inverse Laplace transforms of various terms in (\ref{4.71}):
\begin{equation}\label{4.39}
  \rho_{M,n,m}(t)=\frac{1}{2\pi i}\int_\Gamma M_{n,m}(s)e^{st}ds,
\end{equation}
\begin{equation}\label{4.40}
  \rho_{f,n,m}(t)=\frac{1}{2\pi i}\int_\Gamma\frac{1}{\mu h_n(s,\mu,R)}\Big[\frac{\mu}{s+1}\frac{\partial\xi_{2,n,m}}{\partial r}\Big|_{r=R}-\frac{\partial\phi_{2,n,m}}{\partial r}\Big|_{r=R}\Big]e^{st}ds,
\end{equation}
\begin{equation}\label{4.41}
  E_{1n}(t)=\frac{1}{2\pi i}\int_\Gamma\frac{1}{\mu h_n(s,\mu,R)}e^{st}ds, \ \ \ \ \ \Big(\frac{1}{\mu h_n(s,\mu,R)}=\widehat{E}_{1n}(s)\Big)
\end{equation}
\begin{equation}\label{4.42}
\begin{split}
  E_{2n}(t)&=\frac{1}{2\pi i}\int_\Gamma\frac{P_n(\sqrt{s+1}R)}{h_n(s,\mu,R)}\cdot\frac{1}{(s+1)RP_n(\sqrt{s+1}R)+\frac{n}{R}+\beta}e^{st}ds\\
  &\triangleq\frac{1}{2\pi i}\int_\Gamma\frac{P_n(\sqrt{s+1}R)}{h_n(s,\mu,R)}\cdot\frac{1}{\phi_n(s)}e^{st}ds,
\end{split}
\end{equation}
where $\Gamma: s=J+i\tau$, $-\infty<\tau<\infty$, and $J>\max\{\mbox{Re}\ (\mbox{roots}\ \mbox{of}\ h_n)\}$. Note that if $f_{n,m}^1\equiv f_{n,m}^2\equiv0$ and $b_{n,m}^1\equiv b_{n,m}^2\equiv b_{n,m}^3\equiv0$, then $\rho_{M,n,m}$ is the function $\rho_{n,m}$. The corresponding $w_{n,m}$ and $q_{n,m}$ will be denoted by $w_{M,n,m}$ and $q_{M,n,m}$. We further define $\rho_M=\sum_{n=0}^\infty\sum_{m=-n}^n\rho_{M,n,m}Y_{n,m}$, and similarly define $w_M=\sum_{n=0}^\infty\sum_{m=-n}^nw_{M,n,m}Y_{n,m}$, $q_M=\sum_{n=0}^\infty\sum_{m=-n}^nq_{M,n,m}Y_{n,m}$.

By (\ref{4.71}), we have
\begin{equation}\label{4.43}
\begin{split}
  \widehat{\rho}_{n,m}(s)=&\widehat{\rho}_{M,n,m}(s)+\frac{n+\beta R}{\beta R}\frac{\beta+RP_0(R)}{\beta RP_0(R)}\varepsilon\widehat{\rho}_{f,n,m}(s)\\
  &+\frac{\beta+RP_0(R)}{\beta RP_0(R)}\varepsilon\big(\widehat{b}_{n,m}^1-\frac{n}{R}\widehat{b}_{n,m}^3\big)\widehat{E}_{1n}+\frac{\beta+RP_0(R)}{\beta P_0(R)}\varepsilon \widehat{b}_{n,m}^2\widehat{E}_{2n}\\
  \triangleq&\widehat{\rho}_{M,n,m}(s)+A_1\varepsilon\widehat{\rho}_{f,n,m}(s)+A_2\varepsilon\big(\widehat{b}_{n,m}^1-\frac{n}{R}\widehat{b}_{n,m}^3\big) \widehat{E}_{1n}+A_3\varepsilon \widehat{b}_{n,m}^2\widehat{E}_{2n}.
\end{split}
\end{equation}
Then it follows that
\begin{equation}\label{0.61}
  \rho_{n,m}=\rho_{M,n,m}+A_1\varepsilon\rho_{f,n,m}+A_2\varepsilon \big(b_{n,m}^1-\frac{n}{R}b_{n,m}^3\big)\ast E_{1n}+A_3\varepsilon b_{n,m}^2\ast E_{2n}.
\end{equation}
Furthermore, we derive that
\begin{equation}\label{4.44}
  |\rho_{n,m}-\rho_{M,n,m}|\leq|\varepsilon|\left\{A_1|\rho_{f,n,m}|+A_2\big|\big(b^1_{n,m}-\frac{n}{R}b^3_{n,m}\big)\ast E_{1n}\big|+A_3|b^2_{n,m}\ast E_{2n}|\right\}.
\end{equation}

To gain the estimate $\rho_{n,m}(t)$, we need to estimate the three terms on the right-hand side of (\ref{4.44}), respectively.

We now proceed to estimate the last two terms on the right-hand side of (\ref{4.44}). To begin with, we shall need an improvement of \cite[Lemma 5.5]{HZH2}.
%We now proceed with $n\neq1$. We shall need an improvement of \cite[Lemma 5.5]{HZH2}.
\begin{lem}\label{l10}
There exists $\delta_0>0$, independent of $n$, such that the real parts of the roots of $h_n(s,\mu,R)$ are less than $-\delta_0 n^2$ for all sufficiently large $n$.
\end{lem}
\noindent{\bf Proof.} %Since the proof is similar with \cite[Lemma 5.5]{HZH2}, we just give the difference here. %$j_{n+1/2,m}>(m-1)\pi+\sqrt{\Big(n+\frac{1}{2}\Big)\Big(n+\frac{5}{2}\Big)}$ can be found in \cite[page 2495]{HZH2}.
%\begin{equation*}
%  \widetilde{\phi}_n(s)=(s+1)R+(\frac{n}{R}+\beta)/P_n(\sqrt{s+1}R),
%\end{equation*}
%and
%\begin{equation*}
%  \phi_n(s)=P_n(\sqrt{s+1}R)\widetilde{\phi}_n(s)=(s+1)RP_n(\sqrt{s+1}R)+\frac{n}{R}+\beta,
%\end{equation*}
% Let
%$$\phi_n(s)=(s+1)RP_n(\sqrt{s+1}R)+\frac{n}{R}+\beta.$$
 $j_{n+1/2,m}$ satisfies
$$j_{n+1/2,m}>(m-1)\pi+\sqrt{\Big(n+\frac{1}{2}\Big)\Big(n+\frac{5}{2}\Big)},$$
which can be found in \cite[page 2495]{HZH2}.

Then by the definition of $P_n(\xi)$, we derive
\begin{equation*}
\begin{split}
    \phi_n\Big(-\frac{1}{3R^2}n^2\Big)&=2\Big(-\frac{1}{3R^2}n^2+1\Big)R\sum_{m=1}^{\infty}\frac{1}{(-\frac{1}{3R^2}n^2+1)R^2+j_{n+1/2,m}^2}+\frac{n}{R}+\beta\\
    &>2\Big(-\frac{1}{3R^2}n^2+1\Big)R\\
   &\quad\cdot\sum_{m=1}^{\infty}\frac{1}{(-\frac{1}{3R^2}n^2+1)R^2+\Big[(m-1)\pi+\sqrt{(n+\frac{1}{2})(n+\frac{5}{2})}\Big]^2}+\frac{n}{R}+\beta\\
    &>2\Big(-\frac{1}{3R^2}n^2+1\Big)R\\
    &\quad\cdot\Big\{\frac{1}{(n+\frac{1}{2})(n+\frac{5}{2})-\frac{1}{3}n^2+R^2}+\sum_{m=2}^{\infty}\frac{1}{(m-1)^2\pi^2+\frac{2}{3}n^2}\Big\}+\frac{n}{R}+\beta\\
    &>2\Big(-\frac{1}{3R^2}n^2+1\Big)R\Big\{\frac{1}{\frac{2}{3}n^2}+\sum_{m=1}^{\infty}\frac{1}{m^2\pi^2+\frac{2}{3}n^2}\Big\}+\frac{n}{R}+\beta\\
    &>2\Big(-\frac{1}{3R^2}n^2+1\Big)R\Big\{\frac{1}{\frac{2}{3}n^2}+\frac{1}{\sqrt{\frac{2}{3}}n}\Big\}+\frac{n}{R}+\beta\\
    &>\frac{n-(1+\frac{2}{3}\sqrt{\frac{3}{2}}n)}{R}+\beta\\
    &>\beta
\end{split}
\end{equation*}
for $n$ sufficiently large. Together with $\phi_n(\beta_1+0)=-\infty$, where $\beta_1=-1-(j_{n+1/2,1}/R)^2$ is the pole of the function $P_n(\sqrt{s+1}R)$, we deduce that the first zero $\gamma_1$ lies in the interval $(-1-(j_{n+1/2,1}/R)^2,-n^2/(3R^2))$. Thus it is possible to choose $\delta_0$ so small that $\gamma_1<-\delta_0 n^2$.

The rest of the proof is similar with that of \cite[Lemma 5.5]{HZH2}, so we omit the details here. \hfill $\Box$
%Let
%\begin{equation*}
%\widetilde{\Sigma}_1=\Big\{s=a+ib;\beta_\ell<a<-\delta_0n^2,|b|<1/R\Big\}.
%\end{equation*}
%As in the proof of \cite[Lemma 5.5]{HZH2}, on $\partial\widetilde{\Sigma}_1\backslash\{s=-\delta_0n+ib;|b|\leq1/R\}$, $|\widetilde{\phi}_n(s)|\geq1$. We now proceed to show $|\widetilde{\phi}_n(s)|\geq1$ on $\{s=-\delta_0n^2+ib;|b|\leq1/R\}$. Since
%\begin{equation*}
%\begin{split}
%    |P_n(\sqrt{s+1}R)|&=\Big|\frac{2}{R^2}\sum_{m=1}^{\infty}\frac{1}{s-\beta_m}\Big|\leq\frac{2}{R^2}\sum_{m=1}^{\infty}\frac{1}{|s-\beta_m|}\\
%    &\leq\frac{2}{R^2}\sum_{m=1}^{\infty}\frac{1}{|-\delta_0n^2-\beta_m|}=2\sum_{m=1}^{\infty}\frac{1}{j_{n+1/2,m}^2+R^2-R^2\delta_0n^2}\\
%    &<2\sum_{m=1}^{\infty}\frac{1}{\big[(m-1)\pi+\sqrt{(n+1/2)(n+5/2)}\big]^2+R^2-R^2\delta_0n^2}\\
%    &<2\sum_{m=1}^{\infty}\frac{1}{(m-1)^2\pi^2+(1-R^2\delta_0)n^2}\\
%    &<2\sum_{m=1}^{\infty}\frac{1}{(m-1)^2\pi^2+\frac{1}{2}n^2}=2\{\frac{2}{n^2}+\sum_{m=1}^{\infty}\frac{1}{m^2\pi^2+\frac{1}{2}n^2}\}<\frac{2\sqrt{2}}{n}
%    %&<2\{\frac{2}{n^2}+\frac{\sqrt{2}}{2n}\}<\frac{2\sqrt{2}}{n},
%\end{split}
%\end{equation*}
%we obtain
%\begin{equation*}
%    |\widetilde{\phi}_n(s)|\geq\Big(\frac{n}{R}+\beta\Big)/|P_n(\sqrt{s+1}R)|-|s+1|R>\frac{n}{2\sqrt{2}}\Big(\frac{n}{R}+\beta\Big)-\sqrt{(\delta_0n^2-1)^2+\frac{1}{R^2}}R\geq1
%\end{equation*}
%for $\delta_0$ small enough. The rest of the proof is similar to the proof of \cite[Lemma 5.5]{HZH2}. \hfill $\Box$

\begin{lem}\label{l4.1}
Let $\mu<\mu_\ast(R)$. For $n\neq1$, there exists a small positive number $\delta$, depending only on $\mu$, $R$ such that
\begin{equation}\label{4.45}
  |E_{1n}|=\Big|\frac{1}{2\pi i}\int_\Gamma\frac{e^{st}ds}{\mu h_n(s,\mu,R)}\Big|\leq Ce^{-\delta(n^3+1)t}+C(n+1)^{-3}e^{-\delta(n^2+1)t},
  %|E_{1n}|=|\frac{1}{2\pi i}\int_\Gamma\frac{e^{st}ds}{\mu h_n(s,\mu,R)}|\leq Ce^{-\delta(n+1)t},
\end{equation}
\begin{equation}\label{4.46}
  |E_{2n}|=\Big|\frac{1}{2\pi i}\int_\Gamma\frac{P_n(\sqrt{s+1}R)}{h_n(s,\mu,R)}\frac{e^{st}}{\phi_n(s)}ds\Big|\leq C(n+1)^{-1}e^{-\delta(n^2+1)t}.
  %|E_{2n}|=|\frac{1}{2\pi i}\int_\Gamma\frac{1}{h_n(s,\mu,R)}\frac{e^{st}}{\widetilde{\phi}_n(s)}ds|\leq Ce^{-\delta(n+1)t}.
\end{equation}
\end{lem}
\noindent{\bf Proof.}\
The function $h_n(s,\mu,R)$ can be rewritten as $h_n(s,\mu,R)=c_1(s+c(n)+k_n(s))$, where $c_1k_n(s)=\frac{R^2P_1(R)+1+\beta R}{(s+1)R+(\frac{n}{R}+\beta)/P_n(\sqrt{s+1}R)}-RP_1(R)$ and $c(n)\approx c_2n^3$ as $n\rightarrow\infty$ with $c_1$, $c_2$ positive constants. By Lemma \ref{l10}, $|k_n(s)|\leq \mbox{const}=c_3$ if $\mbox{Re}\ s>-\delta(n^2+1)$, provided $\delta$ is small enough. Then we can use the same argument as in the proof of \cite[Lemma 4.1]{FH2}, and obtain the estimate (\ref{4.45}).

We now proceed to prove (\ref{4.46}). By (\ref{0.8}), we have $|P_n(\sqrt{s+1}R)|\leq\frac{C}{\sqrt{|s|+1}}$ for $\mbox{Re}\ s\geq-\delta(n^2+1)$. Since $|\frac{1}{\phi_n(s)}|\leq C$ if $\mbox{Re}\ s\geq-\delta(n^2+1)$, it follows that
\begin{equation*}
\begin{split}
  |E_{2n}|=\Big|\frac{1}{2\pi i}\int_\Gamma\frac{P_n(\sqrt{s+1}R)}{h_n(s,\mu,R)}\frac{e^{st}}{\phi_n(s)}ds\Big|=&C\Big|\frac{1}{2\pi i}\int_{-\delta(n^2+1)-i\infty}^{-\delta(n^2+1)+i\infty}\frac{P_n(\sqrt{s+1}R)}{s+c(n)+k_n(s)}\frac{e^{st}}{\phi_n(s)}ds\Big|\\
  %\leq&C\frac{1}{2\pi i}\int_{-\delta(n^2+1)-i\infty}^{-\delta(n^2+1)+i\infty}\Big|\frac{P_n(\sqrt{s+1}R)}{h_n(s,\mu,R)}\Big|e^{st}ds\\
  \leq&Ce^{-\delta(n^2+1)t}\int_{-\infty}^\infty\frac{d\tau}{(|\tau|+n^3+1)(|\tau|^{1/2}+n+1)}\\
  \leq& C(n+1)^{-\frac{4}{3}}e^{-\delta(n^2+1)t},
\end{split}
\end{equation*}
which proves (\ref{4.46}).
 \hfill $\Box$

Furthermore, by Lemma \ref{l4.1}, we have the next result.
\begin{lem}\label{l4.2}
If
\begin{equation*}
  \int^\infty_0e^{2\delta_1t}|b(t)|^2dt\leq A, \ \ \ \ 0<\delta_1<\delta,
\end{equation*}
then for all $n\neq1$,
\begin{equation}\label{4.47}
  \int_0^\infty e^{2\delta_1t}|b\ast E_{jn}|^2dt\leq CA(n+1)^{-6}, \ \ \ j=1,2.
\end{equation}
\end{lem}
The proof of this lemma is the same as in \cite[Lemma 4.2]{FH2}.

From Lemma \ref{l4.2} and (\ref{4.21})-(\ref{4.23}), we immediately obtain the following estimates of the last two terms on the right-hand side of (\ref{4.44}).
\begin{lem}\label{l4.20}
For $\delta_1>0$,
\begin{equation}\label{4.48}
  |\varepsilon|(n+1)^7\int_0^\infty e^{2\delta_1t}\big|\big(b^1_{n,m}-\frac{n}{R}b^3_{n,m}\big)\ast E_{1n}\big|^2dt\leq C(B^1_{n,m}+B^3_{n,m}),
\end{equation}
\begin{equation}\label{4.49}
  |\varepsilon|(n+1)^7\int_0^\infty e^{2\delta_1t}|b^2_{n,m}\ast E_{2n}|^2dt\leq CB^2_{n,m},
\end{equation}
where $C$ is independent of $n$.
\end{lem}
%For readers' convenience, we give \cite[Lemma 4.2]{FH2} as follows:
%\begin{lem}\label{l4.2}
%If
%\begin{equation*}
%  \int^\infty_0e^{2\delta_1t}|b(t)|^2dt\leq A, \ \ \ \ 0<\delta_1<\delta,
%\end{equation*}
%then for all $n\neq1$,
%\begin{equation}\label{4.47}
%  \int_0^\infty e^{2\delta_1t}|b\ast E_{jn}|^2dt\leq CA(n+1)^{-2}, \ \ \ j=1,2.
%\end{equation}
%\end{lem}

Next, we shall establish the estimate of $\rho_{f,n,m}$ on the right-hand side of (\ref{4.44}) as follows:
\begin{lem}\label{l4.3}
For $\delta_1>0$,
\begin{equation}\label{4.55}
  |\varepsilon|(n+1)^8\int_0^\infty e^{2\delta_1t}|\rho_{f,n,m}|^2dt\leq C(F_{n,m}^1+F_{n,m}^2),
\end{equation}
where the constant $C$ is independent of $n$.
\end{lem}
\noindent{\bf Proof.}\
We begin with the term
\begin{equation*}
  L_1\equiv\frac{1}{2\pi i}\int_\Gamma\frac{1}{\mu h_n(s,\mu,R)}\frac{\mu}{s+1}\frac{\partial\xi_{2,n,m}}{\partial r}\Big|_{r=R}e^{st}ds,
\end{equation*}
which appears in the definition of $\rho_{f,n,m}$. It follows from (\ref{4.300}) that $\xi_{2,n,m}$ is the Laplace transform of the solution $\Psi_{2,n,m}$ of
\begin{equation*}
  \partial_t\Psi_{2,n,m}-\Delta\Psi_{2,n,m}+\Big(\frac{n(n+1)}{r^2}+1\Big)\Psi_{2,n,m}=f^1_{n,m} \ \ \ \ \mbox{in}\ B_R,\ t>0,
\end{equation*}
\begin{equation*}
  \Big(\frac{\partial\Psi_{2,n,m}}{\partial r}+\beta \Psi_{2,n,m}\Big)\Big|_{r=R}=0, \ \ \ \ \Psi_{2,n,m}|_{t=0}=0.
\end{equation*}
Since
%\begin{equation*}
%\begin{split}
%  \int_0^\infty \mu e^{-t}\ast\Psi_{2,n,m}e^{-st}dt&=\mu\int^\infty_0\int_0^te^{-(t-\tau)}\Psi_{2,n,m}(r,\tau)e^{-st}d\tau dt\\
%  &=\mu\int^\infty_0e^{-(s+1)t}\int_0^te^{\tau}\Psi_{2,n,m}(r,\tau)d\tau dt\\
%  &=\mu\int^\infty_0\int_\tau^\infty e^{-(s+1)t}e^{\tau}\Psi_{2,n,m}(r,\tau)dtd\tau\\
%  &=\frac{\mu}{s+1}\int^\infty_0e^{-s\tau}\Psi_{2,n,m}(r,\tau)d\tau=\frac{\mu}{s+1}\xi_{2,n,m},
%\end{split}
%\end{equation*}
the Laplace transform of $\mu e^{-t}\ast\Psi_{2,n,m}$ is $\frac{\mu}{s+1}\xi_{2,n,m}$, then we can write
\begin{equation*}
  L_1=E_{1n}\ast\frac{\partial}{\partial r}(\mu e^{-t}\ast\Psi_{2,n,m})|_{r=R}.
\end{equation*}
By parabolic estimates and (\ref{4.19}), as in the proof of \cite[(91) of Lemma 5.3]{HZH2}, we obtain
%\begin{equation*}
%  (n+1)^2\Big|\frac{\partial}{\partial r}\Psi_{2,n,m}(R,t)\Big|^2\leq C\|f_{n,m}^1(\cdot,t)\|^2_{L^2(B_R)},
%\end{equation*}
%and then by (\ref{4.19}), we obtain
\begin{equation}\label{4.50}
  |\varepsilon|\int_0^\infty e^{2\delta_1t}(n+1)^2\Big|\frac{\partial\Psi_{2,n,m}(R,t)}{\partial r}\Big|^2dt\leq CF_{n,m}^1.
\end{equation}
Let $b\equiv\mu e^{-t}\ast\frac{\partial\Psi_{2,n,m}}{\partial r}|_{r=R}$, then, by changing the order of integration, we have
\begin{equation*}
\begin{split}
   |\varepsilon|\int_0^\infty e^{2\delta_1t}|b|^2dt=&|\varepsilon|\mu^2\int_0^\infty e^{2\delta_1t}\Big|\int_0^t\frac{\partial\Psi_{2,n,m}(R,t-\tau)}{\partial r}e^{-\tau}d\tau\Big|^2dt\\
   =&|\varepsilon|\mu^2\int_0^\infty e^{2\delta_1t}\Big|\int_0^t\frac{\partial\Psi_{2,n,m}(R,t-\tau)}{\partial r}e^{-\frac{\tau}{2}+\frac{\delta_1\tau}{2}}e^{-\frac{\tau}{2}+\frac{\delta_1\tau}{2}}e^{-\delta_1\tau}d\tau\Big|^2dt\\
   \leq&|\varepsilon|\mu^2\int_0^\infty e^{2\delta_1t}\Big\{\int_0^t\Big|\frac{\partial\Psi_{2,n,m}(R,t-\tau)}{\partial r}\Big|^2e^{-\tau+\delta_1\tau}e^{-2\delta_1\tau}d\tau\int_0^te^{-\tau+\delta_1\tau}d\tau\Big\}dt\\
   \leq&C|\varepsilon|\mu^2\int_0^\infty \int_0^t\Big|\frac{\partial\Psi_{2,n,m}(R,t-\tau)}{\partial r}\Big|^2e^{2\delta_1(t-\tau)}e^{-\tau+\delta_1\tau}d\tau dt\\
   \leq&C|\varepsilon|\mu^2\int_0^\infty \int_\tau^\infty\Big|\frac{\partial\Psi_{2,n,m}(R,t-\tau)}{\partial r}\Big|^2e^{2\delta_1(t-\tau)}e^{-\tau+\delta_1\tau}dtd\tau \\
   =&C|\varepsilon|\mu^2\int_0^\infty \int_0^\infty\Big|\frac{\partial\Psi_{2,n,m}(R,t)}{\partial r}\Big|^2e^{2\delta_1t}e^{-\tau+\delta_1\tau}dtd\tau \\
   \leq& CF^1_{n,m}(n+1)^{-2},
\end{split}
\end{equation*}
where we have used (\ref{4.50}). Hence, it follows from Lemma \ref{l4.2} that
\begin{equation}\label{4.51}
\begin{split}
  |\varepsilon|\int_0^\infty e^{2\delta_1t}|L_1|^2dt=|\varepsilon|\int_0^\infty e^{2\delta_1t}|E_{1n}\ast b|^2dt
  \leq CF^1_{n,m}(n+1)^{-8}.
\end{split}
\end{equation}

Now we consider the term
\begin{equation*}
  L_2\equiv\frac{1}{2\pi i}\int_\Gamma\frac{1}{\mu h_n(s,\mu,R)}\frac{\partial\phi_{2,n,m}}{\partial r}\Big|_{r=R}e^{st}ds.
\end{equation*}
By (\ref{4.36}), it is easily seen that $\phi_{2,n,m}$ is the Laplace transform of the solution $\Phi_{2,n,m}$ of
\begin{equation}\label{4.52}
  -\Delta\Phi_{2,n,m}+\frac{n(n+1)}{r^2}\Phi_{2,n,m}=f^2_{n,m}+\mu e^{-t}\ast f_{n,m}^1 \ \ \ \ \mbox{in}\ B_R,\ t>0,
\end{equation}
\begin{equation}\label{4.53}
  \Big(\frac{\partial\Phi_{2,n,m}}{\partial r}+\beta\Phi_{2,n,m}\Big)\Big|_{r=R}=0,
\end{equation}
then $L_2$ can be rewritten as
\begin{equation*}
  L_2=E_{1n}\ast\frac{\partial\Phi_{2,n,m}(R,t)}{\partial r}.
\end{equation*}
From \cite[Lemma 5.2]{HZH2}, it follows that
\begin{equation*}
  (n+1)^2\Big|\frac{\partial\Phi_{2,n,m}(R,t)}{\partial r}\Big|^2\leq C\big[\|f_{n,m}^2(\cdot,t)\|^2_{L^2(B_R)}+\|e^{-t}\ast f_{n,m}^1(\cdot,t)\|^2_{L^2(B_R)}\big].
\end{equation*}
By the same argument as before, we get
\begin{equation*}
  |\varepsilon|\int_0^\infty e^{2\delta_1t}\|e^{-t}\ast f_{n,m}^1(\cdot,t)\|_{L^2(B_R)}^2dt\leq CF_{n,m}^1,
\end{equation*}
so that
\begin{equation*}
  |\varepsilon|\int_0^\infty e^{2\delta_1t}\Big|\frac{\partial\Phi_{2,n,m}(R,t)}{\partial r}\Big|^2dt\leq C(F_{n,m}^1+F_{n,m}^2)(n+1)^{-2}.
\end{equation*}
Furthermore, it follows from Lemma \ref{l4.2} that
\begin{equation*}\label{4.54}
  |\varepsilon|\int_0^\infty e^{2\delta_1t}|L_2|^2dt=|\varepsilon|\int_0^\infty e^{2\delta_1t}\Big|E_{1n}\ast\frac{\partial\Phi_{2,n,m}(R,t)}{\partial r}\Big|^2dt\leq C(F_{n,m}^1+F_{n,m}^2)(n+1)^{-8}.
\end{equation*}
Combining this estimate with (\ref{4.51}), we derive (\ref{4.55}). Hence, the proof is complete.
\hfill $\Box$

By Lemmas \ref{l4.20} and \ref{l4.3}, we obtain an estimate for (\ref{4.44}) as follows:
\begin{lem}\label{l4.4}
For all $n\neq1$ and $|m|\leq n$:
\begin{equation}\label{4.70}
  (n+1)^7\int_0^\infty e^{2\delta_1t}|\rho_{n,m}(t)-\rho_{M,n,m}(t)|^2dt\leq C|\varepsilon|(F_{n,m}^1+F_{n,m}^2+B_{n,m}^1+B_{n,m}^2+B_{n,m}^3),
\end{equation}
and therefore also
\begin{equation}\label{4.56}
  \int_0^\infty e^{2\delta_1t}\|\rho-\rho_M-\sum_m(\rho_{1,m}-\rho_{M,1,m})Y_{1,m}\|^2_{H^{7/2}(\partial B_R)}dt\leq C|\varepsilon|.
\end{equation}
\end{lem}
As in the proof of \cite{FH2}, we derive the estimates for $w$, $q$ and $\rho_t$:
\begin{equation}\label{4.57}
  \int_0^\infty e^{2\delta_1t}\|w-w_M-\sum_m(w_{1,m}-w_{M,1,m})Y_{1,m}\|^2_{H^{2}( B_R)}dt\leq C|\varepsilon|,
\end{equation}
\begin{equation}\label{4.58}
  \int_0^\infty e^{2\delta_1t}\|q-q_M-\sum_m(q_{1,m}-q_{M,1,m})Y_{1,m}\|^2_{H^{2}( B_R)}dt\leq C|\varepsilon|,
\end{equation}
\begin{equation}\label{4.59}
  \int_0^\infty e^{2\delta_1t}\Big\|\frac{\partial}{\partial t}(\rho-\rho_M-\sum_m(\rho_{1,m}-\rho_{M,1,m})Y_{1,m})\Big\|^2_{H^{1/2}(\partial B_R)}dt\leq C|\varepsilon|.
\end{equation}

\subsection{Case 2: $n=1$}
As stated in Section 4, a translation of the origin does not change the equations (\ref{3.15})-(\ref{3.18}), but change the initial data, so that in the new coordinate system,  by (\ref{0.66}), $\widehat{\rho}_{1,m}(s)$ is changed into
\begin{equation}\label{0.75}
\begin{split}
    \widehat{\rho}_{1,m}(s)
    =&\frac{1+\beta R}{\beta R}\frac{\beta+RP_0(R)}{\beta RP_0(R)}\frac{1}{\mu h_1(s,\mu,R)}\\
    &\cdot\Big\{\frac{\beta R}{1+\beta R}\widetilde{\rho}_{0,1,m}+\frac{\mu}{s+1}\frac{\partial\xi_{1,1,m}}{\partial r}\Big|_{r=R}-\frac{\partial\phi_{1,1,m}}{\partial r}\Big|_{r=R}\\
    &\quad+\varepsilon\Big[\frac{\mu}{s+1}\frac{\partial\xi_{2,1,m}}{\partial r}\Big|_{r=R}-\frac{\partial\phi_{2,1,m}}{\partial r}\Big|_{r=R}+\frac{\beta R}{1+\beta R}\widehat{b}^1_{1,m}\\
  &\quad+\frac{\beta R^2}{1+\beta R}\mu\widehat{b}^2_{1,m}\frac{1}{(s+1)R+(\frac{1}{R}+\beta)/P_1(\sqrt{s+1}R)}-\frac{\beta }{1+\beta R}\widehat{b}^3_{1,m}\Big]\Big\},
\end{split}
\end{equation}
or
\begin{equation}\label{0.76}
\begin{split}
  \widehat{\rho}_{1,m}(s)=&M_{1,m}(s)+\frac{1+\beta R}{\beta R}\frac{\beta+RP_0(R)}{\beta RP_0(R)}\frac{\varepsilon}{\mu h_1(s,\mu,R)}\\
  &\cdot\Big\{\frac{\mu}{s+1}\frac{\partial\xi_{2,1,m}}{\partial r}\Big|_{r=R}-\frac{\partial\phi_{2,1,m}}{\partial r}\Big|_{r=R}+\frac{\beta R}{1+\beta R}\widehat{b}^1_{1,m}\\
  &\quad+\frac{\beta R^2}{1+\beta R}\mu\widehat{b}^2_{1,m}\frac{1}{(s+1)R+(\frac{1}{R}+\beta)/P_1(\sqrt{s+1}R)}-\frac{\beta }{1+\beta R}\widehat{b}^3_{1,m}\Big\},
\end{split}
\end{equation}
where
\begin{equation}\label{0.77}
\begin{split}
  M_{1,m}(s)=&\frac{1+\beta R}{\beta R}\frac{\beta+RP_0(R)}{\beta RP_0(R)}\frac{1}{\mu h_1(s,\mu,R)}\\
  &\cdot\Big\{\frac{\beta R}{1+\beta R}\widetilde{\rho}_{0,1,m}+\frac{\mu}{s+1}\frac{\partial\xi_{1,1,m}}{\partial r}\Big|_{r=R}-\frac{\partial\phi_{1,1,m}}{\partial r}\Big|_{r=R}\Big\},
\end{split}
\end{equation}
and the singularity of $1/h_1(s)$ at $s=0$ is cancelled by the expression in the brace in (\ref{0.75}) at $s=0$, guaranteed by results from section 4. Furthermore, we can take the inverse Laplace transform in (\ref{0.75}) and move the contour $\Gamma$ to $\mbox{Re}\ s=-\delta$, thus obtaining the formula
\begin{equation}\label{5.1}
  \rho_{1,m}(t)=\frac{1}{2\pi i}\int_{-\delta-i\infty}^{-\delta+i\infty}\widehat{\rho}_{1,m}(s)e^{st}ds.
\end{equation}
In order to estimate $\rho_{1,m}(t)$ in (\ref{5.1}), we break up $\widehat{\rho}_{1,m}(s)$ given by the right-hand side of (\ref{0.76}) into four terms, and each term has the form
\begin{equation}\label{5.2}
  \frac{A(s)}{h_1(s)},
\end{equation}
then we shall estimate the inverse Laplace transform separately for each term. However, since $s=0$ is a simple root of $h_1(s)$ for each term (\ref{5.2}), we therefore rewrite (\ref{5.2}) in the form
\begin{equation}\label{5.3}
  \Big(\frac{s}{s+1}\frac{1}{h_1(s)}\Big)A(s)+\Big(\frac{s}{s+1}\frac{1}{h_1(s)}\Big)B(s), \ \ \ B(s)=\frac{A(s)}{s}.
\end{equation}

Introduce new functions
\begin{equation*}
  E_{1\ast}=\frac{1}{2\pi i}\int_{-\delta-i\infty}^{-\delta+i\infty}\frac{s}{s+1}\frac{e^{st}}{\mu h_1(s)}ds,
\end{equation*}
\begin{equation*}
  E_{2\ast}=\frac{1}{2\pi i}\int_{-\delta-i\infty}^{-\delta+i\infty}\frac{s}{s+1}\frac{1}{\mu h_1(s)}\frac{P_1(\sqrt{s+1}R)}{(s+1)RP_1(\sqrt{s+1}R)+\frac{1}{R}+\beta}e^{st}ds.
\end{equation*}
Since the singularity of $1/h_1(s)$ at $s=0$ is cancelled by the factor $s/(s+1)$ at $s=0$ and $s/(s+1)$ goes to $1$ as $|s|\rightarrow\infty$, by the same method as in the proof of Lemma \ref{l4.1} we immediately obtain the following estimates
\begin{equation}\label{5.4}
  |E_{1\ast}|=\Big|\frac{1}{2\pi i}\int_{-\delta-i\infty}^{-\delta+i\infty}\frac{s}{s+1}\frac{e^{st}}{\mu h_1(s)}ds\Big|\leq Ce^{-\delta t},
\end{equation}
\begin{equation}\label{5.5}
  |E_{2\ast}|=\Big|\frac{1}{2\pi i}\int_{-\delta-i\infty}^{-\delta+i\infty}\frac{s}{s+1}\frac{e^{st}}{\mu h_1(s)}\frac{P_1(\sqrt{s+1}R)}{(s+1)RP_1(\sqrt{s+1}R)+\frac{1}{R}+\beta}e^{st}ds\Big|\leq Ce^{-\delta t}.
\end{equation}
Then we can write $\widehat{\rho}_{1,m}(s)$ in the form as
\begin{equation}\label{5.6}
\begin{split}
  \widehat{\rho}_{1,m}(s)=&\widehat{\rho}_{M,1,m}(s)+\frac{1+\beta R}{\beta R}\frac{\beta+RP_0(R)}{\beta RP_0(R)}\varepsilon\widehat{\rho}_{f,1,m}(s)\\
  &+\frac{\beta+RP_0(R)}{\beta R^2P_0(R)}\varepsilon(R\widehat{b}_{1,m}^1(s)-\widehat{b}_{1,m}^3)\widehat{E}_{1\ast}(s)\Big(1+\frac{1}{s}\Big)\\
  &+\frac{\beta+RP_0(R)}{\beta P_0(R)}\varepsilon\widehat{b}_{1,m}^2(s)\widehat{E}_{2\ast}(s)\Big(1+\frac{1}{s}\Big).
\end{split}
\end{equation}
{ For the first term $\widehat{\rho}_{M,1,m}(s)$, we have
%\begin{equation*}
%  \rho_{M,1,m}(t)=\widetilde{\rho}_{M,1,m}(t)+c_m\varepsilon
%\end{equation*}
\begin{equation}
  |\rho_{M,1,m}(t)|=\Big|\frac{1}{2\pi i}\int_{-\delta-i\infty}^{-\delta+i\infty}\widehat{\rho}_{M,1,m}(s)e^{st}ds\Big|=\widetilde{\rho}_{M,1,m}(t)+c_m(\varepsilon),
  %\leq Ce^{-\delta t}.
\end{equation}
where $|c_m(\varepsilon)|\leq C|\varepsilon|$, and by the linear stability result \cite[Theorem 1.2]{HZH2}, $|\widetilde{\rho}_{M,1,m}(t)|\leq Ce^{-\delta t}$ holds.}
%\begin{equation*}
%  |\widetilde{\rho}_{M,1,m}(t)|\leq Ce^{-\delta t}.
%\end{equation*}}

We then apply the previous lemmas and \cite[Lemmas 5.2-5.4]{FH2}  to estimate the inverse Laplace transform of the last three terms in (\ref{5.6}). For example, for the term
\begin{equation*}
\begin{split}
  |\varepsilon|\widehat{b}_{1,m}^2(s)\widehat{E}_{2\ast}(s)\Big(1+\frac{1}{s}\Big)=&|\varepsilon|\widehat{b}_{1,m}^2(s)\widehat{E}_{2\ast}(s)
  +|\varepsilon|\widehat{b}_{1,m}^2(s)\widehat{E}_{2\ast}(s)\frac{1}{s},
  %\equiv& A+B.
\end{split}
\end{equation*}
since the inverse Laplace transform of $\widehat{b}_{1,m}^2\widehat{E}_{2\ast}$ is $b_{1,m}^2\ast E_{2\ast}$, for the first term on the right-hand side we can use the same argument as in the proof of Lemma \ref{l4.2} to establish
\begin{equation}\label{0.70}
  |\varepsilon|\int_0^\infty e^{2\delta_1t}|b_{1,m}^2\ast E_{2\ast}|^2\leq CB_{1,m}^2,
\end{equation}
while for the last term on the right-hand side, by \cite[Lemmas 5.3 and 5.4]{FH2} and (\ref{0.70}), we have
\begin{equation*}
\begin{split}
  |\varepsilon|\int_0^\infty e^{2\delta_1t}\Big|\frac{1}{2\pi i}\int_{-\delta-i\infty}^{-\delta+i\infty}\widehat{b}_{1,m}^2\widehat{E}_{2\ast}\frac{e^{st}}{s}ds\Big|^2dt
  &=|\varepsilon|\int_0^\infty e^{2\delta_1t}\Big|\frac{1}{2\pi i}\int_{-\delta-i\infty}^{-\delta+i\infty}\widehat{b_{1,m}^2\ast E_{2\ast}}\frac{e^{st}}{s}ds\Big|^2dt\\
  &=|\varepsilon|\int_0^\infty e^{2\delta_1t}\Big|\int_t^\infty b_{1,m}^2\ast E_{2\ast}d\tau\Big|^2dt\\
  &\leq\frac{1}{\delta_1^2}\int_0^\infty e^{2\delta_1t}|b_{1,m}^2\ast E_{2\ast}|^2dt\\
  &\leq CB_{1,m}^2.
\end{split}
\end{equation*}

Together with the above analysis, the main result for the case $n=1$ is established.
\begin{lem}\label{l5.1}
For a positive constant $C$, we get
\begin{equation}\label{5.7}
  \int_0^\infty e^{2\delta_1t}|\rho_{1,m}(t)-\rho_{M,1,m}(t)|^2dt\leq C|\varepsilon|.
\end{equation}
In particular,
\begin{equation}\label{5.8}
  \int_0^\infty e^{2\delta_1t}|\rho_{1,m}(t)|^2dt\leq C.
\end{equation}
\end{lem}

\section{Stability for $\mu<\mu^\ast$}
We shall establish the asymptotic stability of the radially symmetric stationary solution for (\ref{3.2})-(\ref{3.6}) by using two fixed point theorem. In Section 4, by the first fixed point theorem (Theorem \ref{th2}) we determine the new center $\varepsilon a^\ast(\varepsilon)$. And as in the proof of \cite[Lemma 7.1]{FH2}, we obtain that after performing a translation $x\rightarrow x+\varepsilon a^\ast(\varepsilon)$ on the initial data, there exists a unique global solution $(w,q,\rho)$ of (\ref{3.15})-(\ref{3.18}), satisfying
\begin{equation}\label{0.30}
\begin{split}
  \|w\|_{C^{2+2\alpha/3,1+\alpha/3}(B_R\times[0,\infty))}\leq C,
\end{split}
\end{equation}
\begin{equation}\label{0.31}
\begin{split}
  \|q\|_{C^{2+\alpha,\alpha/3}(B_R\times[0,\infty))}\leq C,
\end{split}
\end{equation}
\begin{equation}\label{0.32}
\begin{split}
  \|\rho,D_x\rho\|_{C^{3+\alpha,1+\alpha/3}(\partial B_R\times[0,\infty))}\leq C.
\end{split}
\end{equation}

We now proceed with a second fixed point argument.

Introduce the space $X$ of functions $\Phi=(f^1,f^2,b^1,b^2,b^3)$ with norm $\|\Phi\|$ defined by the maximum of the left-hand sides of (\ref{4.1})-(\ref{4.10}) with $\sqrt{|\varepsilon|}$ dropped, and set
\begin{equation*}
  X_1=\{\Phi\in X:\ \sqrt{|\varepsilon|}\|\Phi\|\leq1\},
\end{equation*}
then we define a new function $\widetilde{\Phi}\equiv S\Phi=(\widetilde{f}^1,\widetilde{f}^2,\widetilde{b}^1,\widetilde{b}^2,\widetilde{b}^3)$, where
\begin{equation*}
\begin{split}
  \widetilde{f}^1=-A_\varepsilon^1w+A_\varepsilon w, \ \ \ \ \widetilde{f}^2=-A_\varepsilon w,\\
  \widetilde{b}^1=B_\varepsilon^1,\ \ \ \ \ \widetilde{b}^2=B_\varepsilon^2,\ \ \ \ \ \widetilde{b}^3=B_\varepsilon^3.
\end{split}
\end{equation*}
%\noindent for $a=a^\ast$,
%\begin{equation*}
%  \widetilde{f}^1(r,\theta,\varphi,t;a^\ast)=-A_\varepsilon^1w+A_\varepsilon w, \ \ \ \ \widetilde{f}^2(r,\theta,\varphi,t;a^\ast)=-A_\varepsilon w,
%\end{equation*}
%\begin{equation*}
%  \widetilde{b}^1(\theta,\varphi,t;a^\ast)=B_\varepsilon^1, \ \ \ \widetilde{b}^2(\theta,\varphi,t;a^\ast)=B_\varepsilon^2, \ \ \ \widetilde{b}^3(\theta,\varphi,t;a^\ast)=B_\varepsilon^3;
%\end{equation*}
%\noindent for $a\neq a^\ast$, let
%\begin{equation*}
%  \overline{f}^1(r,\theta,\varphi,t;a)=-A_\varepsilon^1w+A_\varepsilon w, \ \ \ \ \overline{f}^2(r,\theta,\varphi,t;a)=-A_\varepsilon w,
%\end{equation*}
%\begin{equation*}
%  \overline{b}^1(\theta,\varphi,t;a)=B_\varepsilon^1, \ \ \ \overline{b}^2(\theta,\varphi,t;a)=B_\varepsilon^2, \ \ \ \overline{b}^3(\theta,\varphi,t;a)=B_\varepsilon^3,
%\end{equation*}
%and define
%\begin{equation*}
%  \widetilde{f}^i(r,\theta,\varphi,t;a^\ast)=(1-\zeta(t))\widetilde{f}^i(r,\theta,\varphi,t;a^\ast)+\zeta(t)\overline{f}^i(r,\theta,\varphi,t;a),
%\end{equation*}
%\begin{equation*}
%  \widetilde{b}^j(r,\theta,\varphi,t;a^\ast)=(1-\zeta(t))\widetilde{b}^j(r,\theta,\varphi,t;a^\ast)+\zeta(t)\overline{b}^j(r,\theta,\varphi,t;a),
%\end{equation*}
%with the cutoff function $\zeta\in C^\infty[0,\infty)$ such that $\zeta(t)=1$ for $t<1/2$, $\zeta(t)=0$ for $t>1$, and $0\leq\zeta\leq1$ for $1/2\leq t\leq1$.
Again as in the proof of \cite{FH2}, we obtain that $S$ maps $X_1$ into itself and $S$ is a contraction mapping so that $S$ has a unique fixed point.
%Namely,
%\begin{equation}\label{0.59}
%  \int_0^\infty e^{2\delta_1t}\|\rho\|^2_{H^{7/2}(\partial B_R)}dt\leq C|\varepsilon|.
%\end{equation}
%From (7.12) in \cite{FH2},
%\begin{equation*}
%  \|\rho(\cdot,T)\|_{C^{4+\alpha}(B_R)}+\|w(\cdot,T)\|_{C^{2+2\alpha/3}(B_R)}\leq4,
%\end{equation*}
%it follows that the solution and the estimate (\ref{0.59}) can be extended to $[T,2T]$, $[2T,3T]$, $[3T,4T]$, $\cdots$.

By the above argument, we can now state the main result of this paper.
\begin{thm}\label{th4}
Let $\mu<\mu^\ast$. If $|\varepsilon|$ is sufficiently small, then there exists a unique global solution $(\sigma,p,r)$ of the problem (\ref{1.2})-(\ref{1.6}) with the initial data (\ref{1.13}) satisfying (\ref{0.57}), and there exists a new center $\varepsilon a^\ast(\varepsilon)$, where $a^\ast(\varepsilon)$ is a bounded function of $\varepsilon$, such that
\begin{equation*}
  \partial\Omega(t)\rightarrow\big\{|x-\varepsilon a^\ast(\varepsilon)|=R\big\}
\end{equation*}
exponentially fast as $t\rightarrow\infty$.
\end{thm}

\begin{rem}
After the translation of the origin and the Hanzawa transformation, the global solution $(\sigma,p,r)$ in the new variables $(r,\theta,\varphi)$ has the form:
\begin{equation*}\label{3.1}
\begin{split}
    \sigma(r,\theta,\varphi,t)&=\sigma_S(r)+\varepsilon w(r,\theta,\varphi,t),\\
    p(r,\theta,\varphi,t)&=p_S(r)+\varepsilon q(r,\theta,\varphi,t),\\
    \partial\Omega(t):\ r&=R+\varepsilon\rho(\theta,\varphi,t),
\end{split}
\end{equation*}
where $w$, $q$, $\rho$ satisfy (\ref{0.30})-(\ref{0.32}).
\end{rem}

\section{Acknowledgments}
 %We would like to thank the referees for their valuable recommendations.
  The second author is supported by the National Natural Science Foundation of China (No. 11371286) and the Natural Science Basic Research Plan in Shaanxi Province of China (No. 2019JM-165).

\end{document}